\documentclass[twoside,10pt]{article}

\usepackage{a4,amsthm,latexsym,amsmath,amssymb,amsfonts}

\def\bs{\bigskip}
\def\ms{\medskip}
\def\1{\mbox{\it 1 \hskip -7pt  I}}

\def\b{\hbox{\bf b}}
\def\sb{\hbox{\bf{\tiny b}}}

\def\E{{\scriptscriptstyle E}}
\def\T{{\scriptscriptstyle T}}
\def\a{\hbox{\bf a}}
\def\e{\hbox{\bf e}}
\def\se{\hbox{\bf{\tiny e}}}

\def\nn{{\mathbb N}}

\def\no{\noindent}
\def\rr{{\mathbb R}}

\def\bb{{\mathfrak B}}

\def\mm{{\mathcal M}}

\def\di{\displaystyle}

\def\d{\delta}
\let\g=\gamma
\def\limn{\lim_{n\rightarrow\infty}}     
\newtheorem{definition}{{\bf{\small D}{\scriptsize EFINITION}}}
\newtheorem{theorem}{{\bf{\small T}{\scriptsize HEOREM}}}
\newtheorem*{main-theorem}{{\bf{\small M}{\scriptsize AIN THEOREM}}}
\newtheorem{lemma}{{\bf{\small L}{\scriptsize EMMA}}}
\newtheorem{proposition}{{\bf{\small P}{\scriptsize ROPOSITION}}}

\newtheorem{remark}{{\bf{\small R}{\scriptsize EMARK}}}
\newtheorem*{notation}{{\bf{\small N}{\scriptsize OTATIONS}}}

\begin{document}

\title{Projection of Markov measures\\ may be Gibbsian}

\maketitle

\centerline{\scshape J.-R. Chazottes$^{a}$, E. Ugalde$^{b}$}
\medskip
 
{\footnotesize
\centerline{$^{a}$ CPHT-CNRS}
\centerline{{\'E}cole Polytechnique, 91128 Palaiseau Cedex, France}
\centerline{{\tt Email: jeanrene@cpht.polytechnique.fr}} 
\centerline{$^{b}$ IICO - UASLP}
\centerline{A. Obreg{\'o}n 64, 78000 San Luis Potos{\'\i}, SLP, M{\'e}xico}
\centerline{{\tt Email: ugalde@cactus.iico.uaslp.mx}}
}

\medskip

\begin{abstract}
We study the induced measure obtained from a 1-step Markov measure, 
supported by a topological Markov chain, after the mapping of the original
alphabet onto another one.
We give sufficient conditions for the induced measure to be a Gibbs
measure (in the sense of Bowen) when the factor system is again a topological
Markov chain. This amounts to constructing, when it does exist,
the induced potential and proving its H{\"o}lder continuity. This is
achieved through a matrix method. 
We provide examples and counterexamples to illustrate our results.
\end{abstract}

\newpage

\tableofcontents

\newpage

%%%%%%%%%%%%%%%%%%%%%%%%%%%%%%%%%%%%%%%%%%%%%%%%%%
%%%%%%%%%%%%%%%%% NEW SECTION
%%%%%%%%%%%%%%%%%%%%%%%%%%%%%%%%%%%%%%%%%%%%%%%%%%
\section{Introduction}

This paper is concerned with the nature of the ``projection'' of a Markov measure,
supported by a topological Markov chain (TMC for short), obtained by the action of
a factor map mapping the original TMC onto another one.
(We recall the definition of a TMC at the beginning of Section \ref{setup}.)
The resulting measure is not expected to be a Markov measure of any order, that is, the
resulting process has not a finite memory.
The simplest class of measures with infinite memory one could expect is the class of
Bowen-Gibbs measures (BGM's for short).

\no This problem arises naturally in the coding under restrictions of
the kind imposed by forbidding the use of certain blocks. A factor map
(called a code in that context) represents a channel with deterministic noise, that is, one
which looses information in a predictable way \cite{MPW}. Here the
input messages are governed by statistics described by a Markov chain 
and one wants to determine statistics of output messages. 

This problem is also related to the so-called hidden Markov model \cite{rabiner} in
Statistics: this model consists in assuming that the observed data are the image
of a finite-state Markov chain, this image being obtained by ``lumping'' some
of the states of the state space. Our problem can be phrased by saying that we
wish to determine whether a hidden Markov process is distributed according to
a BGM.

A third situation where our problem naturally arises is the following. Suppose that
a chaotic time series $\{x_n\}_{n\geq 0}$ is generated by a deterministic process,
a dynamical system, and assume for the sake of definiteness that it is a map on the interval.
This means that $x_{n+1}=f(x_n)$. In general one does not have access to $f$ and only observes
a symbolic sequence instead of the original orbit. This is because the system can only be
observed through a partition of the values of the $x_n$'s that corresponds to the finite precision
of the measurement or the computer. A natural question is to determine the invariant measure
from this single symbolic sequence, even in the ideal case when the time series would be of
infinite length. This problem has been for instance studied in \cite{cfl} where it was assumed
that the observed symbolic sequence is generating by a Bowen-Gibbs measure. A particular class
of maps $f$ is the one of piecewise linear Markov maps.
When one consider the coding of such maps
via the partition given by the intervals corresponding to each branch, the resulting symbolic
dynamics is given by a TMC with a state space with $k$ symbols and
the invariant measure is a ($1$-step) Markov measure, $k$ being the number of branches \cite{KH}.
A basic question is the following: if one
observes the dynamical system through a lumping of the partition just introduced, supposing
that two atoms of the partition cannot be distinguished, say, then what is the resulting
invariant measure describing the time-series ?

A last incarnation of our problem is a one-dimensional lattice gas described by a Markov
measure. What happens if, say, two spin values cannot be distinguished ? What we call a projection
un the present article is in that context an example of a single site renormalization group transformation.
Non Gibbsianess is not expected since there are no phase transitions in one dimensional finite
range systems. Some useful references for the reader interested in classical models of statistical
mechanics are \cite{enter,maes2,maes3,maes1}. Of course, while we restrict ourselves to the
one-dimensional setting, the problem of transforming Gibbs measures (by many other types
of procedures) can be set in the much more general context of measures on $d$-dimensional lattices,
see \cite{fernandez} for the most recent review. 

\bs 

%%%%%%%%%%%%%%%%%%%%%%%%%%%%%%%%%%%%%%%%%%%%%%%%%%%%%%%%%%%%%%%%
{\bf Description of the paper.} 
Section \ref{setup} is devoted to the set-up of our article. We also give
the ansatz for the induced potential based on a simple property of a 
Gibbs measure. It turns out that the point is to control an infinite
product of non-square matrices.

In section \ref{mainresult} we state our main result, namely some sufficient
conditions to get a BGM from the original 1-step Markov measure after the
projection of its state space. We emphasize that the presence of forbidden
blocks in the original system causes the main difficulty. The projection
process induces some strong topological correlations in the resulting
system and the existence of the ansatz potential is not obvious at all. 

The main result is proved in section \ref{proofofmaintheorem}.
In a first subsection we define a suitable projective metric which
is the central tool to control the infinite products of matrices
appearing in the ansatz of the induced potential. In the following subsection
we state a theorem giving some sufficient conditions on a point in the
projected TMC in order to have a well-defined potential at that point. In the
last subsection we show how to extend the preceding theorem to the whole
projected TMC and we prove the H{\"o}lder continuity of the induced potential.
Therefore, under suitable conditions, the projection of the initial Markov
measure is a BGM.

Section \ref{examples} provides a typical example illustrating our main result.
Then we also consider the case when the original TMC is a full shift, that is
when no blocks are forbidden. It turns out that the projected measure
is always a BGM, generically with an infinite range potential (in very special
cases the potential can be of finite range).
Notice that in the absence of forbidden blocks our problem is considerably simplified.
We also present an example showing that one of the two hypothesis needed to establish our
main result is not just technical. Indeed in that example the induced potential is not
defined at some point (the infinite product mentioned above does not converge). This
also illustrates the non-trivial effect that the presence of forbidden blocks in the
original system may have. Then we give a formula of the induced potential, when it
is well-defined, at periodic points. This follows from the classical Perron-Frobenius
theorem since we have for such points to perform products of positive square matrices.

In section \ref{open} we give some links between our paper and some related works
both in our context and in other settings. We also address some natural issues
raised by our main result and the counterexamples. 

A last section contains the proof of some auxiliary lemmas.

%%%%%%%%%%%%%%%%%%%%%%%%%%%%%%%%%%%%%%%%%%%%%%%%%%%%%%%%%%%%%%%%%%%%
%%%%%%%%%%%%%%%%% SUBSECTION %%%%%%%%%%%%%%%%%%%%%%%%%%%%%%%%%%%%%%%
\section{Set-up and ansatz for the induced potential}\label{setup}

%%%%%%%%%%%%%%%%%%%%%%%%%%%%%%%%%%%%%%%%%%%%%%%%%%%%%%%%%%%%%%%%%%%%
\subsection{Set-up}

Let $(A_M,\sigma)$ be a TMC where $A$ stands for the (finite) alphabet, $M$ for the incidence 
matrix and $\sigma$ denotes the shift map. This means that $M$ is a $0-1$-matrix selecting
a subset of all possible infinite sequences $\a=\a(0)\a(1)\cdots$ drawn from
the alphabet $A$:
$$
A_M :=\{\a\in A^{\nn}: M(\a(i),\a(i+1))=1\;\forall\;i\in\nn\}\, .
$$
This subset is closed under the action of the shift transformation $\sigma$,
that is, $\sigma A_M=A_M$, where $\sigma$ is defined as follows:
$(\sigma\a)(i):=\a(i+1)$ for any $\a\in A_M$. (TMC's are nothing but
subshifts of finite type with forbidden block of length two.)
Let us recall that a TMC can be viewed as the set of infinite paths
on the directed graph (digraph) with vertex set equal to the
alphabet and arrows corresponding to allowed transitions between symbols
of the alphabet according to the incidence matrix \cite{K}. (We will use this
representation in section \ref{examples}.)

Suppose $(A_M,\sigma)$ is {\bf topologically mixing} or, equivalently, that $M$ is
a {\bf primitive matrix}.
By this we mean that there is a power $m_0\geq 1$ such that the matrix
$M^{m_0}$ has only strictly positive entries. 
(Notice that this property is equivalent to assume that $M$ is irreducible and
aperiodic \cite{seneta}.)
Consider a 1-step Markov measure $\mu:\bb(A_M)\to [0,1]$, which 
is $\sigma$--invariant and mixing ($\bb(A_M)$ denotes the Borel 
sigma-algebra of $A_M$ generated by cylinder sets). 
This measure is a BGM associated to a potential $\phi:A_M\to \rr$ which is constant 
in each cylinder of length two. This potential can be thought as a 2-symbols function, that
is a potential of range two.
Using the same notation for both the potential and the 2-symbols function,
we have $\phi(\a)=\phi(\a(0:1))$, where $\a(0:1)$ denotes the 2-block $\a(0)\a(1)$.
In general, given $\a\in A_M$ and $0\leq i < j$, the block $\a(i)\a(i+1)\cdots
\a(j)$ will be denoted by $\a(i:j)$.

\no The potential generating $\mu$ can always be chosen such that
%%%%%%%%%%%%%%%%%%%%%%%%%%%%%%%%%%%%%%%%%%%%%%%%%%%%%%%%
\begin{equation}\label{markov-potential}
\mu[\a(0:n)]=\exp\left(\sum_{i=0}^{n-1}\phi(\a(i:i+1))\right)\mu[\a(n)]\, ,
\end{equation}
%%%%%%%%%%%%%%%%%%%%%%%%%%%%%%%%%%%%%%%%%%%%%%%%%%%%%%%%
where $[\a(0:n)]$ is the cylinder of length $(n+1)$ containing $\a$, i.~e., 
$[\a(0:n)]:=\{\a'\in A_M:\ \a'(i)=\a(i),\ 0\leq i\leq n\}$. (In this case
$\phi$ must be strictly negative.)

\bs {\bf Bowen-Gibbs inequality}. Let us recall the following basic characterization
of a general BGM. Let $\Omega$ be a TMC and $\psi: \Omega\to \rr$. Then it is
known \cite{bowen} that there is a unique $\sigma$-invariant measure $\nu$ such
that for any $n\in\nn_0$ and any admissible $\b$
%%%%%%%%%%%%%%%%%%%%%%%%%%%%%%%%%%%%%%%%%%%%%%%%%%%%%%%%%%%%%%%%%%%%%%%
\begin{equation}\label{BGI}
\exp(-K)\leq
\frac{\nu[\b(0:n)]}{\exp\left(\sum_{j=0}^{n}\psi(\sigma^j(\b))\right)}
\leq  \exp(K). 
\end{equation}
%%%%%%%%%%%%%%%%%%%%%%%%%%%%%%%%%%%%%%%%%%%%%%%%%%%%%%%%%%%%%%%%%%%%%%%
where $K>0$ is a constant independent of $n$ and $\b$. We implicitely
put the topological pressure of $\psi$ equal to zero, which is always
possible \cite{pp}.

\bs \no From now on, we always choose inside the cohomological class
of potentials determining any BGM the normalized one, so, in particular, the
one of zero pressure \cite{pp}.

\bs \no Let $\pi:A\to B$ be a map onto another alphabet $B$. In the sequel
we always assume that $\#A>\#B>1$, that is $\pi$ is a projection. 
This map defines a factor system $(B_\pi,\sigma)$, with
%%%%%%%%%%%%%%%%%%%%%%%%%%%%%%%%%%%%%%%%%%%%%%%%%%%%%%%
\begin{equation}\label{factor-system}
B_{\pi}:=\{\b\in B^\nn:\ \exists\ \a\in A_M \hbox{ such that } \pi
\a(i)=\b(i),\ i\in\nn_0\}\, .
\end{equation}
%%%%%%%%%%%%%%%%%%%%%%%%%%%%%%%%%%%%%%%%%%%%%%%%%%%%%%%%
($\nn_0= \nn\cup \{0\}$. We also denote by $\pi$ the map from $A_M$ onto $B_\pi$ defined in an obvious way.)
  
\bs \no It is readily checked that if $B_\pi$ is infinite then it
has to be uncountable. Notice that $(B_\pi,\sigma)$ is in general a
topologically mixing sofic subshift \cite{K}. A sofic subshift cannot
be described by a list of finite forbidden blocks.

\bs \no Let us introduce the following distance on $B_\pi$
%%%%%%%%%%%%%%%%%%%%%%%%%%%%%%%%%%%%%%%%%%%%%%%%%%%%%%%%%
\begin{equation}
\label{defmetric}
d(\b,\b'):=\left\{
\begin{array}{lcr}
{\di 
\exp\left(-\frac{\min\{j:\b(j)\neq\b'(j)\}}{2(\# B+1)}\right) } & 
                                       \textup{if} & \b\neq \b' \\
                                      0 & \textup{if} & \b =\b'.
\end{array}\right.
\end{equation}
%%%%%%%%%%%%%%%%%%%%%%%%%%%%%%%%%%%%%%%%%%%%%%%%%%%%%%
(The constant dividing the min in this definition is for the
sake of later convenience.)

\bs \no The problem concerning us is to elucidate the nature of 
the measure $\nu:\bb(B_\pi) \to [0,1]$ such that
%%%%%%%%%%%%%%%%%%%%%%%%%%%%%%%%%%%%%%%%%%%%%%%%%%%%%%%%%
\begin{equation}\label{induced-measure}
\nu[\b(0:n)]:=\mu\{\a\in A_M:\ \pi(\a(i))=\b(i),\ 
\forall \ 0\leq i\leq n\},
\end{equation}
%%%%%%%%%%%%%%%%%%%%%%%%%%%%%%%%%%%%%%%%%%%%%%%%%%%%%%%%
which is the image (or projection) of the measure $\mu$ by $\pi$, i.e. 
$\nu:=\mu\circ \pi^{-1}$.

\bs\no Some general properties of the original measure are preserved under the action of factor maps.
If $\mu$ is ergodic so is the measure $\nu\circ\pi^{-1}$. The same holds for
the mixing property. (See e.g. \cite{KH}.)

\bs\no For each $b\in B$ let $E_{b}:=\pi^{-1}(b)\subset A$.
For each $bb'\in B\times B$ define the rectangular matrix
$\mm_{bb'}:E_{b}\times E_{b'}\to [0,1]$ by
%%%%%%%%%%%%%%%%%%%%%%%%%%%%%%%%%%%%%%%%%%%%%%%%%%%%%%%%%%%%%%
\begin{equation}\label{b(0:1)-submatrices}
\mm_{bb'} (a,a')=\exp\left[\phi(a,a')\right]M(a,a')
\end{equation}
%%%%%%%%%%%%%%%%%%%%%%%%%%%%%%%%%%%%%%%%%%%%%%%%%%%%%%%%%%%%%%
for all  $aa'\in E_{b}\times E_{b'}$.

\bs \no Finally, for each $b\in B$ let us define the column vector
%%%%%%%%%%%%%%%%%%%%%%%%%%%%%%%%%%%%%%%%%%%%%%%%%%%%%%%%%%%%%%
\begin{equation}\label{mub}
\mu_b:\, E_b\to [0,1]\quad\text{such that}\quad\mu_b(a)=\mu[a],\,
\text{for any}\, a\in E_b\,.
\end{equation}
%%%%%%%%%%%%%%%%%%%%%%%%%%%%%%%%%%%%%%%%%%%%%%%%%%%%%%%%%%%%%%

\bs\no A straightforward computation shows that
$\nu:\bb(B_\pi)\to[0,1]$ satisfies
%%%%%%%%%%%%%%%%%%%%%%%%%%%%%%%%%%%%%%%%%%%%%%
\begin{equation}\label{induced-via-matrices}
\nu[b_0 b_1\cdots b_n]=\1^{\dag}\left(\prod_{i=0}^{n-1}
\mm_{b_i b_{i+1}}\right)\mu_{b_n}
\end{equation}
%%%%%%%%%%%%%%%%%%%%%%%%%%%%%%%%%%%%%%%%%%%%%%
for any $B_\pi$--admissible cylinder $[b_0 b_1\cdots b_n]$.
Here the symbol $\1$ stands for the all-ones column vector of the 
adequate dimension (whereas $\1^{\dag}$ is the corresponding all-ones 
row vector).

%%%%%%%%%%%%%%%%%%%%%%%%%%%%%%%%%%%%%%%%%%%%%%%%%%%%%%%%%%%%%%%%%%%%%%
\subsection{Ansatz for the induced potential}

\no In the following basic lemmas $\Omega$ denotes an arbitrary TMC. 

\bs
%%%%%%%%%%%%%%%%%%%%%%%%%%%%%%%%%%%%%%%%%%%%%%%%%%%%%%%%%%%%%%%%
\begin{lemma}\label{potential-lemma}
Suppose that $\nu:\bb(\Omega)\to [0,1]$ is a BGM. Then it is
associated to the normalized H{\"o}lder continuous potential
$\psi:\Omega\to\rr$
such that
%%%%%%%%%%%%%%%%%%%%%%%%%%%%%%%%%%%%%%%%%%%%%%%%%%%%%%
\begin{equation}\label{markovapproximation}
\psi(\b)=\lim_{n\to\infty}\log\left(\frac{\nu[\b(0:n)]}
{\nu[\b(1:n)]}\right)
\end{equation}
%%%%%%%%%%%%%%%%%%%%%%%%%%%%%%%%%%%%%%%%%%%%%%%%%%%%%%
for any $\b\in \Omega$.
\end{lemma}
%%%%%%%%%%%%%%%%%%%%%%%%%%%%%%%%%%%%%%%%%%%%%%%%%%%%%%%%%%%%

\no We refer to~\cite{pp} for the straightforward proof of this result.
We remark that
$\psi_n(\b):=\log(\nu[\b(0:n)]/\nu[\b(1:n)])$
defines for each $n\geq 1$ a function which is constant in each cylinder
of length $n+1$. This function is the H{\"o}lder continuous potential of a certain
(unique) BGM $\nu_n$, say, which is nothing but the $n$-step Markov approximation of the
measure $\nu$.

\bs\no Let us emphasize that Lemma \ref{potential-lemma} says that {\it necessarily}
the potential of a BGM is given by \eqref{markovapproximation}. This is far from sufficient
since this lemma holds for any $g$-measure. A $g$-measure is an equilibrium state (which is
in general not unique) associated to a suitably normalized stictly positive potential which 
is only continuous. We refer the reader to \cite{ppw} and references therein for details.

\bs \no {\bf Ansatz for the induced potential}. Coming back to our problem, we see that 
Lemma \ref{potential-lemma} gives the
following ansatz for the induced potential: for each $\b\in B_\pi$, set

%%%%%%%%%%%%%%%%%%%%%%%%%%%%%%%%%%%%%%%%%%%%%%%%%%%%%%
\begin{eqnarray}\label{the-limit}
\nonumber
\psi(\b) & = & \lim_{n\to\infty}\log
\left(\frac{\nu[\b(0:n)]}{\nu[\b(1:n)]}\right)\\
         & = &\lim_{n\to\infty}
\log\left( \frac{\1^{\dag}\left(\prod_{i=0}^{n-1}
\mm_{\sb(i:i+1)}\right)\mu_{\sb(n)}
    }{\1^{\dag}\left(\prod_{i=1}^{n-1}\mm_{\sb(i:i+1)}\right)\mu_{\sb(n)}}\, ,
\right)     
\end{eqnarray}
%%%%%%%%%%%%%%%%%%%%%%%%%%%%%%%%%%%%%%%%%%%%%%%%%%%%%%%%%%%%%%%%%%%%%%%

\no where the second equality is obtained after the straightforward substitutions
according to formula \eqref{induced-via-matrices}. We see that this ansatz potential
is given by an {\em infinite product of non square matrices} whose convergence and
regularity properties as a function of $\b$ seem to be not trivial at all.
The key tool to control this infinite product will be the use of a suitable projective metric. 

\bs \no In the next section we state (and prove in section \ref{proofofmaintheorem}) that,
under some sufficient conditions, the measure $\nu$ defined by \eqref{induced-measure} is the BGM
associated to the above ansatz potential. This means that we will show that $\psi(\b)$ is correctly
defined for every $\b\in B_{\pi}$ and, moreover, that it is a H{\"o}lder continuous function.
 
%%%%%%%%%%%%%%%%%%%%%%%%%%%%%%%%%%%%%%%%%%%%%%%%%%%%%%%%%%%%%%%%%%
%%%%%%%%%%%%%%% NEW SECTION
%%%%%%%%%%%%%%%%%%%%%%%%%%%%%%%%%%%%%%%%%%%%%%%%%%%%%%%%%%%%%%%%%%
\section{Main result: when a Markov measure is mapped to a BGM}\label{mainresult}

\no The next theorem gives sufficient conditions ensuring that 
the function $\b\mapsto \psi(\b)$, defined by (\ref{the-limit}), is well defined and H{\"o}lder continuous 
in the whole projected TMC, $B_\pi$. Before stating the theorem, we need some preliminary
definitions.

\bs 
%%%%%%%%%%%%%%%%%%%%%%%%%%%%%%%%%%%%%%%%%%%%%%%%%%%%%%
\begin{definition}[Row allowable matrix]\label{row-allow}
Let $E'$ and $E$ be finite alphabets, and 
$T:E'\times E\to [0,\infty)$ be a rectangular, non--negative matrix 
on these alphabets. This matrix is said to be row allowable if 
for each $e'\in E'$ there exists $e\in E$ such that $T(e',e) > 0$.
\end{definition}
%%%%%%%%%%%%%%%%%%%%%%%%%%%%%%%%%%%%%%%%%%%%%%%%%%%%%%

\no This definition is inspired by a very similar one given 
in~\cite{seneta}.

\bs

\no First, we restrict the type of factor maps or projections because we want 
to get a TMC from the original one. Otherwise, as mentioned above, one would get
a sofic subshift in which ``topological correlations'' are generally ``non-local''.
  
%%%%%%%%%%%%%%%%%%%%%%%%%%%%%%%%%%%%%%%%%%%%%%%%%%%%%%%%%%%%%%%%%%%%%%%%%%%
\begin{definition}[Topological Markov factor map] The factor map 
$\pi:A_M\to B_\pi$ is said to be a topological Markov map if the factor 
subshift $B_\pi$ is a TMC.
\end{definition}
%%%%%%%%%%%%%%%%%%%%%%%%%%%%%%%%%%%%%%%%%%%%%%%%%%%%%%%%%%%%%%%%%%%%%%%%%%%

\bs \no For each $b\in B$ recall that
\begin{equation} 
E_b:=\{a\in A:\ \pi(a)=b\}\ . \nonumber
\end{equation}

\no For each $B_\pi$--admissible block $bb'$, let 
$M_{bb'}:E_b\times E_{b'}\to \{0,1\}$ be such that
\begin{equation}
M_{bb'}(a,a')=\left\{ \begin{array}{lr}
      1 & \ \text{if} \ M(a,a')=1, \ \pi(a)=b\ \text{and}\ \pi(a')=b',\\
      0 & \ \text{otherwise}. \end{array}\right.
\nonumber
\end{equation}

\no We need a further restriction on the factor maps we will be able to handle.

%%%%%%%%%%%%%%%%%%%%%%%%%%%%%%%%%%%%%%%%%%%%%%%%%%%%%%%%%%%%%%%%%%
\begin{definition}[Full row allowable factor map]
The factor map $\pi:A_M\to B_\pi$ is said to be full row allowable if 
for each $B_\pi$--admissible block $bb'$, the corresponding transition
submatrix $M_{bb'}$ is row allowable. 
\end{definition}
%%%%%%%%%%%%%%%%%%%%%%%%%%%%%%%%%%%%%%%%%%%%%%%%%%%%%%%%%%%%%%%%%%%%%%

\no Notice that $M_{bb'}$ is a submatrix of the transition matrix $M$, 
which is compatible with the non-negative matrix $\mm_{bb'}$ defined 
by (\ref{b(0:1)-submatrices}). It is clear that the definition of a full
row allowable factor map does not depend on the potential $\phi$ defining
the Markov measure on $(A_M,\sigma)$. It is a purely topological notion.

\no Denote by ${\rm Per}_p(B_\pi)$ the set of admissible periodic points
with period $p\geq 1$. Notice that there is at least one $p$ with $1\leq p\leq \#B$ 
such that ${\rm Per}_p(B_\pi)\neq\emptyset$. This is a basic property of a TMC \cite{K}.

\bs\no We can now state the main theorem of the paper.
%%%%%%%%%%%%%%%%%%%%%%%%%%%%%%%%%%%%%%%%%%%%%%%%%%%%%%%%%%%%%%%%%%

\begin{main-theorem}
\label{holder-continuity}
Suppose that
$(A_M,\sigma)$ is a TMC supporting a 1-step
Markov measure $\mu$.
Let $\pi:A\to B$ a map from $A$ onto
another alphabet $B$ and $B_\pi$ the corresponding
factor space. Assume that $\#A>\#B>1$ and
$\pi$ satisfies the following conditions:
\begin{itemize}
\item[{\bf (H1)}] is full row allowable, 
\item[{\bf (H2)}] for each $\b\in {\rm Per}_p(B_\pi)$, with 
$1\leq p \leq \# B$, the matrix
%$\mm_{\sb(0:p)}$
$\prod_{i=0}^{p-1}\mm_{\sb(i:i+1)}$ is positive. 
\end{itemize}

\no Under these hypothesis, the function $\b\mapsto \psi(\b)$
defined by \eqref{the-limit} is well-defined and H{\"o}lder continuous
on the whole set $B_\pi$ which is a TMC. This amounts to saying that the projected measure
$\nu=\mu\circ \pi^{-1}$ (remember formula (\ref{induced-measure}))
is the (unique) BGM of the potential $\psi$, that is,
it satisfies the Bowen-Gibbs inequality \eqref{BGI}. 
\end{main-theorem}
%%%%%%%%%%%%%%%%%%%%%%%%%%%%%%%%%%%%%%%%%%%%%%%%%%%%%%%%%%%%%%%%%%

\no In section \ref{examples} we provide typical examples of factor maps satisfying
the hypotheses of this theorem. Moreover we give a formula of the induced potential
at periodic points. It is natural to ask what happens in the case when $A_M$ is
a full shift. We shall show in that section that the projected measure is always a BGM because
the hypotheses of our theorem are trivially fulfilled in that case. A more interesting
question is whether the corresponding potential can be of finite range (which gives a
Markov measure with a certain memory). A typical 
example will show this is possible but ``non generic''.

\no It is worth to point out
that the presence of forbidden blocks makes the induced potential of infinite range.
The reason is that a potential is of finite range if and only if the sequence in
formula \eqref{markovapproximation} becomes constant after some $n_0$, which means
that the potential is of range $n_{0}+1$ (see the lines just after Lemma \ref{potential-lemma}).
The presence of forbidden blocks makes
unlikely this phenomenon to occur (see formula \eqref{the-limit}).

\no An example built in section \ref{examples} will show that the hypothesis {\bf H2} is
unavoidable since we will exhibit a point such that the induced potential does not exist.

\begin{remark}
We have only considered the case when the original system is a TMC, instead
of a more general subshift of finite type (SFT),
and the factor map is only a 1-block factor map, instead of, say, a 2-block factor
map. A SFT is a subshift for which is given a list of forbidden blocks whose length
is two in the case of TMC's.
From the mathematical point of view there is no loss of generality since any
SFT can be recoded as a TMC and a finite-block factor map as a 1-block factor map.
We are not able to handle the case of a generic sofic subshift.
We refer to \cite{K} for background informations on symbolic
dynamics and coding.
\end{remark}

%%%%%%%%%%%%%%%%%%%%%%%%%%%%%%%%%%%%%%%%%%%%%%%%%%%%%%%%%%%%%%%%%%%%%%%%%%%%%%%%%%%%%%%%%%%%%%%%%%
%%%%%%%%%%%%%%%%%%%%%%%%%%%% NEW SECTION %%%%%%%%%%%%%%%%%%%%%%%%%%%%%%%%%%%%%%%%%%%%%%%%%%%%%%%%%
%%%%%%%%%%%%%%%%%%%%%%%%%%%%%%%%%%%%%%%%%%%%%%%%%%%%%%%%%%%%%%%%%%%%%%%%%%%%%%%%%%%%%%%%%%%%%%%%%%
\section{Proof of the main theorem}\label{proofofmaintheorem}

This section is divided into three subsections. We first introduce the projective
metric we need to control the infinite product of non-square matrices
that appears in \eqref{the-limit}. This is the crucial point in our approach.
Then we give sufficient conditions for a $\b\in B_\pi$ such that $\psi(\b)$
defined in \eqref{the-limit} does exist (Theorem \ref{existence-limit}). 
Then we prove our main theorem.

%%%%%%%%%%%%%%%%%%%%%%%%%%%%%%%%%%%%%%%%%%%%%%%%%%%%%%%%%%%
\subsection{Contractivity of positive non-square matrices over simplices}

\no Let $E\subset A$ be non--empty, and consider the simplex
%%%%%%%%%%%%%%%%%%%%%%%%%%%%%%%%%%%%%%%%%%%%%%%%%%%%%%%%%%
\begin{equation}\label{the-simplex}
\Delta_{\E}:=\left\{x:E\to (0,1): \ \vert x\vert_1=1 \right\}\, ,
\end{equation}
%%%%%%%%%%%%%%%%%%%%%%%%%%%%%%%%%%%%%%%%%%%%%%%%%%%%%%%%%%%%%%
\no where $\vert x\vert_1:=\1^{\dag} x = \sum_{e\in E}x(e)$.

\bs \no The projective metric in this simplex is the function
$\d_{\E}:\Delta_{\E}\times \Delta_{\E}\to [0,\infty)$ such that
%%%%%%%%%%%%%%%%%%%%%%%%%%%%%%%%%%%%%%%%%%%%
\begin{equation}\label{projective-distance}
\d_{\E}(x,y) := \log \left(
\frac{\di \max_{e\in E} (x(e)/y(e))}{\di \min_{e\in E} (x(e)/y(e))}
\right)\,\cdot
\end{equation}
%%%%%%%%%%%%%%%%%%%%%%%%%%%%%%%%%%%%%%%%%%%%

\no The projective metric makes 
$(\Delta_E, \d_{\E})$ a complete metric space.
Of course $\Delta_{\E}$ is not 
complete with respect to the Euclidean, or any other $\ell_p$ metric.

\bs\no Let us associate to any matrix $T:E'\times E\to [0,\infty)$
the mapping $F\!\!_{\T}:\Delta_{\E}\to \Delta_{\E'}$, such that
%%%%%%%%%%%%%%%%%%%%%%%%%%%%%%%%%%%%%%%%%%%%%%%%%%%%
\begin{equation}\label{mapping-on-simplex}
F\!\!_{\T}(x):=\frac{ Tx }{\vert Tx\vert_1 }\ \cdot
\end{equation}
%%%%%%%%%%%%%%%%%%%%%%%%%%%%%%%%%%%%%%%%%%%%%%%%%%%%
Notice that this mapping is well defined for any row allowable
matrix.

\bs 
%%%%%%%%%%%%%%%%%%%%%%%%%%%%%%%%%%%%%%%%%%%%%%%%%%%%
\begin{lemma}\label{contraction}
Let $E, E'\subset A$ be non--empty and
$T:E'\times E\to [0,\infty)$ be row allowable, then
\begin{equation}\d_{\E'}(F\!\!_{\T}(x),F\!\!_{\T}(y))\leq \d_{\E}(x,y),
\end{equation}
\no for all $x,y\in \Delta_{E}$.
\end{lemma}
%%%%%%%%%%%%%%%%%%%%%%%%%%%%%%%%%%%%%%%%%%%%%%%%%%%%%

\bs
%%%%%%%%%%%%%%%%%%%%%%%%%%%%%%%%%%%%%%%%%%%%%%%%%%%%%%%%%
\begin{proof}
We follow the standard procedure one can find in~\cite{seneta}.

\bs \no For $y\in \Delta_{\E'}$ let $P_y:E'\times E\to [0,1]$ be 
such that
%%%%%%%%%%%%%%%%%%%%%%%%%%%%%%%%%%%%%%%%%%%%%%%%%%%%%%%
\begin{equation}
P_y (e',e):=\frac{T(e',e)y(e)}{(Ty)(e')}\ \cdot \nonumber
\end{equation}
%%%%%%%%%%%%%%%%%%%%%%%%%%%%%%%%%%%%%%%%%%%%%%%%%%%%%%%%%
\no for any $(e',e)\in E'\times E$.

\bs \no Then, for all $e'\in E'$ we have
%%%%%%%%%%%%%%%%%%%%%%%%%%%%%%%%%%%%%%%%%%%%%%%%%%%%%%%%
\begin{equation}
\frac{(F\!\!_{\T}(x))(e')}{(F\!\!_{\T}(y))(e')}=
\frac{\vert Ty\vert_1}{\vert Tx\vert_1}
\times \sum_{e\in E} \left(\frac{P_y(e',e)x(e)}{y(e)}\right) =
\left(P_y \left(\frac{x}{y}\right)\right)(e'), \nonumber
\end{equation}
%%%%%%%%%%%%%%%%%%%%%%%%%%%%%%%%%%%%%%%%%%%%%%%%%%%%%%
\no where ${\di \left(\frac{x}{y}\right)}$ denotes the vector 
of quotients ${\di \frac{y(e)}{x(e)}\in (0,\infty)^E}$. 

\no Since $P_y$ is a stochastic matrix, then
\begin{eqnarray}
\max_{e'\in E'}\frac{(Tx)(e')}{(Ty)(e')}&\leq &
                    \max_{e\in E}\frac{x(e)}{y(e)}, \nonumber \\
              \min_{e'\in E'}\frac{(Tx)(e')}{(Ty)(e')}&\geq &
                    \min_{e\in E}\frac{x(e)}{y(e)}, \nonumber
                    \end{eqnarray}
\no which implies that 
$\d_{\E'}(F\!\!_{\T}(x),F\!\!_{\T}(y))\leq \d_{\E}(x,y)$.
\end{proof}
%%%%%%%%%%%%%%%%%%%%%%%%%%%%%%%%%%%%%%%%%%%%%%%%%%%%%%%%%

\no According to this lemma, a non--negative matrix $T$ defines a
``non--expanding'' map $:F\!\!_{\T}:\Delta_E\to \Delta_{E'}$
for which one defines a ``contraction coefficient''.

%%%%%%%%%%%%%%%%%%%%%%%%%%%%%%%%%%%%%%%%%%%%%%%%%%%%%%%%%%%%%%%%%%%%%%%%%%%%
\begin{definition}\label{def-contraction-coefficient}
Let $:F\!\!_{\T}:\Delta_E\to \Delta_{E'}$ be
the mapping defined above (formula \eqref{mapping-on-simplex}).
Then this mapping is ``contractive'' with a
``contraction coefficient'' equal to 
%%%%%%%%%%%%%%%%%%%%%%%%%%%%%%%%%%%%%%%%%%%%%%%%%%%%%%%
\begin{equation}
\label{formula-contracion-coefficient}
\tau(T):=\sup_{x,y\in
\Delta_{\E}}\frac{\d_{\E'}(F\!\!_{\T}(x),F\!\!_{\T}(y))}{\d_{\E}(x,y)}
\cdot 
\end{equation}
%%%%%%%%%%%%%%%%%%%%%%%%%%%%%%%%%%%%%%%%%%%%%%%%%%%%%%%
\end{definition}
%%%%%%%%%%%%%%%%%%%%%%%%%%%%%%%%%%%%%%%%%%%%%%%%%%%%%%%%%%%%%%%%%%%%%%%%%%%%

\no According to Lemma \ref{contraction}, this coefficient is never larger
than one and under certain conditions it is strictly smaller. In
fact one can derive an expression for the contraction coefficient.

%%%%%%%%%%%%%%%%%%%%%%%%%%%%%%%%%%%%%%%%%%%%%%%%%%%%%%%%%%
\begin{lemma}[Contraction coefficient]\label{contraction-coefficient}
For $T:E'\times E\to [0,\infty )$ as above, 
the contration coefficient \eqref{formula-contracion-coefficient}
is equal to
\begin{equation}
\tau(T)=\frac{1-\sqrt{\Phi(T)}}{1+\sqrt{\Phi(T)}}
\nonumber
\end{equation}

\no with
\begin{equation}
\nonumber
\Phi(T)=\left\{\begin{array}{cr}
 {\di \min_{e,f\in E \atop  e',f'\in E'}
\frac{T(e',e)T(f',f)}{T(e',f)T(f',e)}}  &  \ {\rm if }\ T>0 \\
                                      0 &  \ {\rm if }\ T\geq 0
\,.\end{array}\right.
\end{equation}
\end{lemma}
%%%%%%%%%%%%%%%%%%%%%%%%%%%%%%%%%%%%%%%%%%%%%%%%%

\bs \no We do not give the proof of this result since
it can be deduced straightforwardly {\it mutatis mutandis}
from~\cite[section 3.4]{seneta} (therein this formula is
deduced in the case of square matrices). Let us stress that effective
contraction is possible only if the matrix is strictly positive.

%%%%%%%%%%%%%%%%%%%%%%%%%%%%%%%%%%%%%%%%%%%%%%%%%%%%%%%%%%
\subsection{Existence of the induced potential at some points}\label{existence}

\no The aim of this section is to determine sufficient conditions on $\b$
under which the limit (\ref{the-limit}) does exist.

\bs \no Before investigating the existence of this limit in the
whole symbolic space $B_\pi$, we shall consider conditions that
ensure its existence for particular choices of $\b\in B_\pi$. 
But before we need to introduce some shorthand notations in order to
avoid cumbersome formulas.

\bs
\begin{notation}
For $\b\in B_\pi$ and for any integers $0\leq m<n$, let
%%%%%%%%%%%%%%%%%%%%%%%%%%%%%%%%%%%%%%%%%%%%%%%%%%%%%%%%%%%%
\begin{equation}\label{Mmn}
\mm_{\sb(m:n)}:=\prod_{i=m}^{n-1}\mm_{\sb(i:i+1)}\,\cdot
\end{equation}
%%%%%%%%%%%%%%%%%%%%%%%%%%%%%%%%%%%%%%%%%%%%%%%%%%%%%%%%%%%%

\bs \no 
Notice that, since $\mu_{\sb(n)} > 0$ (remember (\ref{mub}))
and $\mm_{\sb(m:n)}\neq 0$, then $\mm_{\sb(m:n)}\mu_{\sb(n)}\neq 0$.

\no We will assume that  $\mm_{\sb(m:n)}$ 
is a row allowable matrix (see Definition \ref{row-allow}).
In this case we have $\mm_{\sb(m:n)}\mu_{\sb(n)}>0$ and
the transformation 
$F_{\mm_{\sb(m:n)}}$ defined by the matrix $\mm_{\sb(m:n)}$
will be simply denoted by $F_{\sb(m:n)}$. 

\no Let us write  
$
x_{\sb(m:n)}:=F_{\sb(m:n)}\left(\hat{\mu}_{\sb(n)}\right)
$
for the image by $F_{\sb(m:n)}$ of the  normalized vector 
$\hat{\mu}_{\sb(n)}:=\mu_{\sb(n)}/|\mu_{\sb(n)}|_1$. 
For $\mm_{\sb(m:n)}$ row allowable, this vector 
lies in the simplex $\Delta_{E_{\sb(m)}}$ (i.~e., 
$x_{\sb(m:n)} > 0$). 

\no Let us denote by $\Delta_{\sb(m)}$ the simplex $\Delta_{E_{\sb(m)}}$, and
by $\d_{\sb(n)}(.,.)$ the corresponding projective distance $\d_{E_{\sb(n)}}(.,.)$.
\end{notation}

%%%%%%%%%%%%%%%%%%%%%%%%%%%%%%%%%%%%%%%%%%%%%%%%%%%%%%%%%%%%%
\bs \no With the notations just introduced we have
%%%%%%%%%%%%%%%%%%%%%%%%%%%%%%%%%%%%%%%%%%%%%%%%%%%%%%%%%%%%
\[
\frac{\1^{\dag}\left(\prod_{i=0}^{n-1}\mm_{\sb(i:i+1)}\right)\mu_{\sb(n)}
}{\1^{\dag}\left(\prod_{i=1}^{n-1}\mm_{\sb(i:i+1)}\right)\mu_{\sb(n)}} = 
\frac{\1^{\dag} \mm_{\sb(0:n)}\mu_{\sb(n)}}{\1^{\dag} \mm_{\sb(1:n)} \mu_{\sb(n)}}=
\]
\[
=
\1^{\dag} \mm_{\sb(0:1)} F_{\sb(1:n)} (\hat{\mu}_{\sb(n)}) 
= \1^{\dag} \mm_{\sb(0:1)}x_{\sb(1:n)}.
\]
%%%%%%%%%%%%%%%%%%%%%%%%%%%%%%%%%%%%%%%%%%%%%%%%%%%%%%%%%%%%%%

\bs \no Therefore, proving that limit~(\ref{the-limit}) exists amounts
to proving that 
\begin{equation}\label{psiasx}
\psi(\b):=\limn \log\left(\1^{\dag} \mm_{\sb(0:1)}x_{\sb(1:n)}\right)
\end{equation}
exists.
In fact, under suitable conditions, $x_{\sb(1:n)}$ converges
exponentially fast to a certain vector, as $n\to\infty$.

\bs \no
%%%%%%%%%%%%%%%%%%%%%%%%%%%%%%%%%%%%%%%%%%%%%%%%%%%%%%%%%%%
\begin{theorem}\label{existence-limit}
Let $\b\in B_\pi$ be such that:
\begin{itemize}
\item[{\bf (H1')}] 
for each $i\geq 0$, 
$\mm_{\sb(i:i+1)}:E_{\sb(i)}\times E_{\sb(i+1)}\to (0,\infty)$ is row allowable; 
\item[{\bf (H2')}] there exists a bounded gap, strictly increasing sequence  
$\ell(0)<\ell(1)<\cdots $ 
(i.~e., $0< \ell(k+1)-\ell(k)<s$ for all $k\in\nn$ and 
some fixed $s\geq 2$), 
such that $\mm_{\sb(\ell(k):\ell(k+1))} > 0$ for each $k \in \nn_0$.
\end{itemize}

\ms \no Under the above conditions, there exists a (normalized) vector
$x_{\sb(1:\infty)}\in \Delta_{\sb(1)}$ and constants $\theta(\b)\in (0,1)$,
$C(\b)>0$, such that 
%%%%%%%%%%%%%%%%%%%%%%%%%%%%%%%%%%%%%%%%%%%%%%%%%%%%%%%%%%%%%
\begin{equation}
\delta_{\sb(1)}(x_{\sb(1:\infty)}, x_{\sb(1:n)})\leq C(\b) \ 
\theta(\b)^{n} \quad \forall n\in \nn.
\end{equation}
%%%%%%%%%%%%%%%%%%%%%%%%%%%%%%%%%%%%%%%%%%%%%%%%%%%%%%%%%%%%%
Moreover,
the induced potential $\psi$ at point $\b\in B_\pi$ exists and it satisfies
%%%%%%%%%%%%%%%%%%%%%%%%%%%%%%%%%%%%%%%%%%%%%%%%%%%%%%%%%%%%%%%
\begin{equation}\label{potential-as-limit}
\left|\psi(\b)-\log\left(\1^{\dag}\mm_{\sb(0:1)}x_{\sb(1:n)}
\right)\right|
      \leq C(\b) \ \theta(\b)^{n} \quad 
\forall n\in \nn.
\end{equation}
\end{theorem}
%%%%%%%%%%%%%%%%%%%%%%%%%%%%%%%%%%%%%%%%%%%%%%%%%%%%%%%%%%%%%%%

\bs\no {\it Proof of Theorem \ref{existence-limit}.}
Hypothesis {\bf H1} implies that each one of the matrices
$\mm_{\sb(m:n)}$ is row allowable (a product of row allowable matrices
being a row allowable matrix).
Therefore, the associated transformation
$F_{\sb(m:n)}:\Delta_{\sb(m)}\to \Delta_{\sb(n)}$ is either an isometry
or a contraction with respect to the suitable projective metric.

\bs\no Hypothesis {\bf H2} and Lemma~\ref{contraction-coefficient} 
imply that $F_{\sb(\ell(k):\ell(k+1))}:\Delta_{\sb(\ell(k+1))}\to
\Delta_{\sb(\ell(k))}$ is a contraction for each $k\in \nn_0$. 
Let us denote  the contraction coefficient
of this transformation, $\tau(\mm_{\sb(\ell(k):\ell(k+1))})$, by $\tau(\b, k)$.
(Recall the definition of the contraction coefficient,
Definition \ref{def-contraction-coefficient}.)

\bs\no
Given $n\in\nn$, let $k(n):=\max\{k\in\nn_0 :\ \ell(k+1)\leq n\}$. 
For any  $n' > n$ one has
%%%%%%%%%%%%%%%%%%%%%%%%%%%%%%%%%%%%%%%%%%%%%%%%%%%%%%%%%%%%%%%%%%%%%%
\begin{eqnarray}\label{tau(k)-inequality}
\d_{\sb(1)}\left(x_{\sb(1:n)},x_{\sb(1,n')}\right)
&=&
\d_{\sb(1)}\left(F_{\sb(1:n)}(\hat{\mu}_{\sb(n)}), 
F_{\sb(1:n)}(x_{\sb(1,n')}) \right) \nonumber\\
&\leq&
\left(\prod_{j=0}^{k(n)}\tau(\b,j)\right)
\times \d_{\sb(n)}\left(\hat{\mu}_{\sb(n)},x_{\sb(n:n')}\right)\,\cdot
\end{eqnarray}
%%%%%%%%%%%%%%%%%%%%%%%%%%%%%%%%%%%%%%%%%%%%%%%%%%%%%%%%%%%%%%%%%%%
(Remember that $\hat{\mu}_{\sb(n)}:=\mu_{\sb(n)}/|\mu_{\sb(n)}|_1$.)

\bs \no The number of different positive matrices
$\mm_{\sb(\ell(k):\ell(k+1))}$, which is at most the number of different
blocks $\b(\ell(i):\ell(i+1))$ occurring in $\b$, is finite. Indeed,
because of the bounded gap condition, this number is not larger than
$(\# B)^{s+1}$, where $s$ is the maximum gap length between consecutive
elements in the sequence $\ell(0)<\ell(1)<\cdots$.
Thus, $\tau(\b):=\sup\{\tau(\b,k):\ k\in \nn_0\}$
is a number strictly smaller than 1. Using the definition
of $k(n)$ and the bounded gap condition we deduce that 
$k(n)\geq \frac{n-\ell(0)}{s}-1$. Hence we can
write~(\ref{tau(k)-inequality}) as follows, 
%%%%%%%%%%%%%%%%%%%%%%%%%%%%%%%%%%%%%%%%%%%%%%%%%%%%%%%%%%%%%%%%%%%%
\begin{eqnarray}
\d_{\sb(1)}\left(x_{\sb(1:n)},x_{\sb(1,n')}\right)
&\leq&
\left(\tau(\b)\right)^{k(n)}\times
\d_{\sb(n)}\left(\hat{\mu}_{\sb(n)},x_{\sb(n:n')}\right)\nonumber\\
&\leq&
C_1(\b)\times \theta(\b)^{n}\times
\d_{\sb(n)}\left(\hat{\mu}_{\sb(n)},x_{\sb(n:n')}\right),\nonumber
\end{eqnarray}
where $\theta(\b):=\tau(\b)^{1/s}$, and 
$C_1(\b):=1/\tau(\b)^{2+\ell(0)/s}$.
%%%%%%%%%%%%%%%%%%%%%%%%%%%%%%%%%%%%%%%%%%%%%%%%%%%%%%%%%%%%%%%%%%%%

\bs\no Hence, $x_{\sb(1:n)}$ converges exponentially fast to a limit 
$x_{\sb(1:\infty)}\in \Delta_{\sb(1)}$, provided
that $\d_{\sb(n)}\left(\hat{\mu}_{\sb(n)},x_{\sb(n:n')}\right)$ 
is bounded with respect to $n$ and $n'>n$.
For this note that
%%%%%%%%%%%%%%%%%%%%%%%%%%%%%%%%%%%%%%%%%%%%%%%%%%%%%%%%%%%%%%%%
\begin{eqnarray}
\d_{\sb(n)}\left(\hat{\mu}_{\sb(n)},x_{\sb(n:n')}\right)
 &  \leq  & \d_{\sb(n)}\left(\hat{\mu}_{\sb(n)},x_{\sb(n:\ell(k(n)+2))}\right)
                                                       \nonumber \\
 &  +  & \d_{\sb(n)}\left(x_{\sb(n:\ell(k(n)+2))},x_{\sb(n:\ell(k(n)+3))}\right) 
                                                        \nonumber \\
 &  +  & \d_{\sb(n)}\left(x_{\sb(n:\ell(k(n)+3))},x_{\sb(n:\ell(k(n)+4))}\right) 
                                                        \nonumber \\
 & \vdots &                                             \nonumber \\
 &  +  & \d_{\sb(n)}\left(x_{\sb(n:\ell(k(n)+p))},x_{\sb(n:n')}\right),   
                                                         \nonumber
\end{eqnarray}
%%%%%%%%%%%%%%%%%%%%%%%%%%%%%%%%%%%%%%%%%%%%%%%%%%%%%%%%%%%%
\no with $p=\max\{k\in \nn_0:\ \ell(k(n)+p)\leq n' \}$.

\bs \no Since each one of the transformations $F_{\sb(\ell(k):\ell(k)+1)}$ 
is  contractive with coefficient $\tau(\b,k)\leq \tau(\b)$, one obtains
%%%%%%%%%%%%%%%%%%%%%%%%%%%%%%%%%%%%%%%%%%%%%%%%%%%%%%%%%%
%\begin{eqnarray}
%\d_{\sb(n)}\left(\hat{\mu}_{\sb(n)},x_{\sb(n:n')}\right)
%                &  \leq  & 
%\d_{\sb(n)}\left(\hat{\mu}_{\sb(n)},x_{\sb(n:\ell(k(n)+2))}\right)
%                  \nonumber \\
% & + & \tau(\b)\
%\d_{\sb(\ell(k(n)+2))}\left(\hat{\mu}_{\sb(\ell(k(n)+2))},
%                               x_{\sb(\ell(k(n)+2):\ell(k(n)+2))}\right)
%                   \nonumber \\
% & + & \tau(\b)^2\ 
%\d_{\sb(\ell(k(n)+3))}\left(\hat{\mu}_{\sb(\ell(k(n)+3))},
%                               x_{\sb(\ell(k(n)+3):m\ell(k(n)+4))}\right)
%                   \nonumber \\
% & \vdots &        \nonumber \\
% & + & \tau(\b)^{p-1}\ 
%\d_{\sb(\ell(k(n)+p))}\left(\hat{\mu}_{\sb(\ell(k(n)+p))},
%                                x_{\sb(\ell(k(n)+p):n')}\right).
%\nonumber
%\end{eqnarray}
%%%%%%%%%%%%%%%%%%%%%%%%%%%%%%%%%%%%%%%%%%%%%%%%%%%%%%%%%%%%%%

%%%%%%%%%%%%%%%%%%%%%%%%%%%%%%%%%%%%%%%%%%%%%%%%%%%%%%%%%%
$$
\d_{\sb(n)}\left(\hat{\mu}_{\sb(n)},x_{\sb(n:n')}\right)\leq
\d_{\sb(n)}\left(\hat{\mu}_{\sb(n)},x_{\sb(n:\ell(k(n)+2))}\right)+
$$
$$
+\tau(\b)\ \d_{\sb(\ell(k(n)+2))}\left(\hat{\mu}_{\sb(\ell(k(n)+2))},
x_{\sb(\ell(k(n)+2):\ell(k(n)+2))}\right)+
$$
$$
+\tau(\b)^2\ 
\d_{\sb(\ell(k(n)+3))}\left(\hat{\mu}_{\sb(\ell(k(n)+3))},
                               x_{\sb(\ell(k(n)+3):m\ell(k(n)+4))}\right)+
$$
$$
\vdots
$$
$$
+\tau(\b)^{p-1}\ 
\d_{\sb(\ell(k(n)+p))}\left(\hat{\mu}_{\sb(\ell(k(n)+p))},
                                x_{\sb(\ell(k(n)+p):n')}\right).
$$
%%%%%%%%%%%%%%%%%%%%%%%%%%%%%%%%%%%%%%%%%%%%%%%%%%%%%%%%%%%%%%

\no Finally, since
%%%%%%%%%%%%%%%%%%%%%%%%%%%%%%%%%%%%%%%%%%%%%%%%%%%%%%%%%%%%%%%
\begin{equation}
D:=\max\left\{
\d_{\sb_0}\left(\hat{\mu}_{\sb_0},x_{\sb_0 \sb_1\cdots \sb_m}\right):
\ 1\leq m < s \ \text{and}\  \b_0 \b_1\cdots \b_m \ \text{is} \ B_{\pi}\!-
\text{admissible}\right\}
\nonumber
\end{equation}
%%%%%%%%%%%%%%%%%%%%%%%%%%%%%%%%%%%%%%%%%%%%%%%%%%%%%%%%%%%%%%%%

\no is finite, one gets
%%%%%%%%%%%%%%%%%%%%%%%%%%%%%%%%%%%%%%%%%%%%%%%%%%%%%%%%%%%%
\begin{equation}
\d_{\sb(1)}\left(x_{\sb(1:n)},x_{\sb(1:n')}\right)\leq
C(\b)\!\times\!\theta(\b)^n, \ \forall n' > n,
\label{contraction-formulae}
\end{equation} 

\no where $C(\b):=\frac{D\times C_1(\sb)}{1-\tau(\sb)}\cdot$. We are
done with the proof of Theorem \ref{existence-limit}.

%%%%%%%%%%%%%%%%%%%%%%%%%%%%%%%%%%%%%%%%%%%%%%%%%%%%%%%%%%%%%%%%%

\bs \no Now turn to prove \eqref{potential-as-limit}.
By definition of the projective
distance and using inequality~(\ref{contraction-formulae}), we obtain
%%%%%%%%%%%%%%%%%%%%%%%%%%%%%%%%%%%%%%%%%%%%%%%%%%%%%%%%%%%%
\begin{eqnarray}
| \underbrace{\log\left(\1^{\dag}\mm_{\sb(0:1)}x_{\sb(1:n')}\right)}_{(*)}-
       \log\left(\1^{\dag}\mm_{\sb(0:1)}x_{\sb(1:n)}\right) 
| &\leq&
\d_{\sb(1)}\left(x_{\sb(1:n)},x_{\sb(1,n')}\right)  \nonumber \\
        &\leq& \left(\frac{D\!\times\! C_1(\b)}{1-\tau(\b)}\right)\!
\times\!\theta(\b)^n \nonumber\\
        &:=& C(\b) \times\!\theta(\b)^n.\nonumber
%\label{convergence-of-potential}
\end{eqnarray}
Hence, $\{ \log\left(\1^{\dag}\mm_{\sb(0:1)}x_{\sb(1:n)}\right)
\}_{n=1}^\infty $ is a Cauchy sequence in $\rr$ converging
exponentially fast to $\psi(\b)$. Formula \eqref{potential-as-limit} follows by taking the
limit $n'\to\infty$ in the term (*) in the previous chain of inequalities.
The proof of Theorem \ref{existence-limit} is now finished.

%%%%%%%%%%%%%%%%%%%%%%%%%%%%%%%%%%%%%%%%%%%%%%%%%%%%%%%%%%%%%%%%%%%%%%
\subsection{Proof of the main theorem}

\no The proof of the Main theorem consists in obtaining a uniform version of formula \eqref{potential-as-limit}
in Theorem \ref{existence-limit}. Before doing this, we need a few more lemmas.

\no First, we state the following lemma.
 
%%%%%%%%%%%%%%%%%%%%%%%%%%%%%%%%%%%%%%%%%%%%%%%%%%%%%%%%%
\begin{lemma}\label{factorizing-b} 
Assume that $B_\pi$ is a TMC. For each $\b\in B_\pi$ there exists a sequence
%%%%%%%%%%%%%%%%%%%%%%%%%%%%%%%%%%%%%%%%%%%%%%%%%%%%%%%%%%%%
\begin{equation}
0\leq m(0)< \ell(0)< m(1) <\ell(1)<\cdots
\nonumber
\end{equation}
%%%%%%%%%%%%%%%%%%%%%%%%%%%%%%%%%%%%%%%%%%%%%%%%%%%%%%%%%%%%%
\no such that for each $k\in\nn$, $\b(m(k))=\b(\ell(k))$ and 
$k(\# B+1)\leq \ell(k) < (k+1)(\# B+1)$.
\end{lemma}
%%%%%%%%%%%%%%%%%%%%%%%%%%%%%%%%%%%%%%%%%%%%%%%%%%%%%%%%%%

\no We need also the following lemma which is a partial converse to
Lemma \ref{potential-lemma}:

\begin{lemma}\label{gibbsian-lemma}
Let $\nu:\bb(\Omega)\to [0,1]$ be a $\sigma$-invariant 
measure. Suppose that there is a summable sequence $\{c(n)\in\rr^+\}_{n=0}^\infty$,
and a H{\"o}lder continuous function $\psi:\Omega\to\rr$ 
such that 
%%%%%%%%%%%%%%%%%%%%%%%%%%%%%%%%%%%%%%%%%%%%%%%%%%%%%%%%
\begin{equation}
\left|\psi(\b)-\log
\left(\frac{\nu[\b(0:n)]}{\nu[\b(1:n)]}\right)\right|\leq
c(n) \nonumber
\end{equation}
%%%%%%%%%%%%%%%%%%%%%%%%%%%%%%%%%%%%%%%%%%%%%%%%%%%%%%%%
Then $\nu$ is the BGM of the (normalized) potential $\psi$.
\end{lemma}
%%%%%%%%%%%%%%%%%%%%%%%%%%%%%%%%%%%%%%%%%%%%%%%%%%%%%%%%%
\no (See section \ref{some-lemmas} for the proof.)

\bs 

\no A last lemma that we need:

\bs
%%%%%%%%%%%%%%%%%%%%%%%%%%%%%%%%%%%%%%%%%%%%%%%%%%%%%%%%%%%%%%%%%%%%%%
\begin{lemma}\label{row->markov}
If $\pi:A_M\to B_\pi$ is full row allowable, then it is a topological
Markov factor map, that is $B_\pi$ is a (primitive) TMC. 
\end{lemma}
%%%%%%%%%%%%%%%%%%%%%%%%%%%%%%%%%%%%%%%%%%%%%%%%%%%%%%%%%%%%%
\no (See section \ref{some-lemmas} for the proof.)

\bs

\no It follows from a simple counterexample presented in
Section~\ref{examples}, that the converse to Lemma \ref{row->markov}
is false.

\bs

\no We are ready for the proof of the main theorem:

\no For $\b\in B_\pi$, let $m(0)<\ell(0)<m(1)<\ell(1)<\cdots$ 
be the sequence whose existence is guaranted by Lemma~\ref{factorizing-b}.  
By hypothesis {\bf H1} and Lemma~\ref{row->markov},
we deduce that there is a periodic point inside each one
of the cylinders $[\b(m(k):\ell(k))]$. Indeed, since the block
$\b(m(k):\ell(k))$ is $B_\pi$-admissible and $(B_\pi, \sigma)$ is a TMC, 
the concatenation $\b(m(k):\ell(k)-1)\b(m(k):\ell(k)-1)\cdots $ is
$B_\pi$-admissible. Thus, by hypothesis {\bf H2} the matrices
$\mm_{\sb(m(k):\ell(k))}$ are all positive.

\no Now by {\bf H1} the matrices $\mm_{\sb(\ell(k-1):m(k))}$ 
are row allowable. It is readily checked that the product of a row 
allowable matrix by a positive matrix is again a positive matrix, which
implies that each one of the matrices $\mm_{\sb(\ell(k):\ell(k+1))}$ is 
positive. 

\bs\no Therefore, for each $\b\in B_\pi$, assumptions {\bf H1'} and {\bf H2'} 
of Theorem~\ref{existence-limit} hold, with a gap constant $s=2(\#B+1)$.
Hence, the induced potential exists in the whole space $B_\pi$. 

\bs\no Now let us establish the H{\"o}lder continuity of $\psi$.
%%%%%%%%%%%%%%%%%%%%%%%%%%%%%%%%%%%%%%%%%%%%%%%%%%%%%%%%%%%%%%%%%

\no Since there are finitely many periodic orbits of period less
than or equal to $\#B$, then
%%%%%%%%%%%%%%%%%%%%%%%%%%%%%%%%%%%%%%%%%%%%%%%%%%%%%%%%%%%%%%%%%%% 
\begin{equation}
\tau:=\sup\left\{\tau(\mm_{\sb(0:p)}): \b \in
{\rm Per}_p(B_\pi)\ \text{for some} \ 1\leq p\leq \#B \right\} < 1.
\nonumber
\end{equation}
%%%%%%%%%%%%%%%%%%%%%%%%%%%%%%%%%%%%%%%%%%%%%%%%%%%%%%%%%%%%%%%%%%

%%%%%%%%%%%%%%%%%%%%%%%%%%%%%%%%%%%%%%%%%%%%%%%%%%%%%%%%%%%%
\bs \no Following the proof of Theorem \ref{existence-limit}, set
$\theta:=\tau^{1/2(\#B+1)}$, 
%%%%%%%%%%%%%%%%%%%%%%%%%%%%%%%%%%%%%%%%%%%%%%%%%%%%%%%%%%%%%%%%%%%%%%%%%%%%%
\begin{equation}
D:=\max\left\{
\d_{b_0}\left(\hat{\mu}_{b_0},x_{b_0 b_1\cdots b_m}\right):
1\leq m < 2(\#B+1) \ b_0b_1\cdots b_m \ \text{is} \ B_{\pi}\!-\!
\text{admissible}\right\},
\nonumber
\end{equation}
%%%%%%%%%%%%%%%%%%%%%%%%%%%%%%%%%%%%%%%%%%%%%%%%%%%%%%%%%%%%%%%%%%%%%%%%%%%%%

\no and $C_1:=1/\tau^{3}$.

\no Inequality~\eqref{contraction-formulae} holds uniformly in $B_\pi$, and we have 
%%%%%%%%%%%%%%%%%%%%%%%%%%%%%%%%%%%%%%%%%%%%%%%%%%%%%%%
\begin{equation}\label{basic-inequality}
\d_{\sb(1)}\left(x_{\sb(1:n)},x_{\sb(1:\infty)}\right)\leq
\left(\frac{D\!\times\! C_1}{1-\tau}\right)\!\times\!\theta^n, 
\nonumber
\end{equation} 
%%%%%%%%%%%%%%%%%%%%%%%%%%%%%%%%%%%%%%%%%%%%%%%%%%%%%%%

\no for all $\b\in B_\pi$ and all $n\in\nn$. Hence, formula \eqref{potential-as-limit} in
Theorem \ref{existence-limit} applies uniformly in $\b$, implying that $\psi(\b)$ exists for all
$\b$ and it is such that
\begin{equation}\label{summable-bound}
\left|\psi(\b)-\log
\left(\frac{\nu[\b(0:n)]}{\nu[\b(1:n)]}\right)\right|\leq
\left(\frac{D\!\times\! C_1}{1-\tau}\right)\!\times\!\theta^n.
\end{equation} 
On the other hand, from inequality~(\ref{basic-inequality}) it also 
follows that 
%%%%%%%%%%%%%%%%%%%%%%%%%%%%%%%%%%%%%%%%%%%%%%%%%%%%%%%%
\begin{equation}
\d_{\sb(1)}\left(x_{\sb'(1:\infty)},x_{\sb(1:\infty)}\right)\leq
\left(\frac{2\!\times\!D\!\times\!C_1}{1-\tau}\right)\!\times\!\theta^n,
\end{equation}
for each $\b\in B_\pi$, each $n\in\nn$, and all $\b' \in [\b(0:n)]$.

%%%%%%%%%%%%%%%%%%%%%%%%%%%%%%%%%%%%%%%%%%%%%%%%%%%%%%%

\no Therefore, because of the definition of the projective distance, 
we obtain
%%%%%%%%%%%%%%%%%%%%%%%%%%%%%%%%%%%%%%%%%%%%%%%%%%%%%%%%%%%%%%
\begin{equation}
|\psi(\b')-\psi(\b)|\leq C\!\times\!\theta^{n}, 
\nonumber
\end{equation}
%%%%%%%%%%%%%%%%%%%%%%%%%%%%%%%%%%%%%%%%%%%%%%%%%%%%%%%%%%%%%

\no with $C:=2\times D\times C_1/(1-\theta)$, which implies that
%%%%%%%%%%%%%%%%%%%%%%%%%%%%%%%%%%%%%%%%%%%%%%%%%%%%%%%%%%%%%%
\begin{equation}\label{psi-is-holder}
\textup{var}_n\psi:= \sup\{\vert \psi(\b)-\psi(\b')\vert: \b_j = \b'_j,
0\leq j\leq n\}\leq C\!\times\!\theta^{n}\ .
\end{equation}

\no According to the definition of the metric~\eqref{defmetric}, 
this means the function $\psi$ is H{\"o}lder continuous with
a H{\"o}lder exponent equal to $\log (1/\tau)$, i.e. the logarithm 
of the inverse of the uniform contraction coefficient.
In this way we prove the existence of a H{\"o}lder continuous potential, 
which does satisfy the hypotheses of Lemma~\ref{gibbsian-lemma} 
because of inequality~(\ref{summable-bound}), and the theorem follows.

%%%%%%%%%%%%%%%%%%%%%%%%%%%%%%%%%%%%%%%%%%%%%%%%%%%%%%%%%%%%%%%%%%%%%%%%%%%%%%%%%%%%%%%%%%%%%%%%
%%%%%%%%%%%%%%%%%% NEW SECTION %%%%%%%%%%%%%%%%%%%%%%%%%%%%%%%%%%%%%%%%%%%%%%%%%%%%%%%%%%%%%%%%%
%%%%%%%%%%%%%%%%%%%%%%%%%%%%%%%%%%%%%%%%%%%%%%%%%%%%%%%%%%%%%%%%%%%%%%%%%%%%%%%%%%%%%%%%%%%%%%%%
\section{Examples, counterexamples and the induced potential at periodic points}\label{examples}

\no In this section we provide some examples illustrating the properties of 
the factor map considered above. 

\bs
\subsection{{\em Ad hoc} example}
Define a class of factor maps by the following construction.

\bs \no {\bf 1.-} Let $B_N$ be a topological Markov chain with
transition matrix $N:B\times B\to \{0,1\}$, and let $A$ be a set
such that $\#A \geq \#B>1$.

\no {\bf 2.-} To each $b\in B$ we associate a non--empty subset
$E_b\subset A$.

\no {\bf 3.-} For each two--block $bb'$ admissible in $B_N$  and for 
each $a\in E_b$, choose a non--empty subset $D_{a,b'}\subset E_{b'}$.

\bs \no A $n$-circuit in $B_N$ (corresponding to a circuit in the
digraph defined by $N$) is a $B_N$--admissible block 
$b_0 b_1\cdots b_n$ such that $b_0=b_n$, and such that $b_i\neq b_j$ 
for all $0\leq i < j < n$. 
Note that there is a finite number of $n$-circuits in $B_N$ for a given $n$.

\bs\no The subset $D_{a,b'}$ in item {\bf 3} has to be non--empty, but it is arbitrary except for 
at least one two--block inside each circuit.

\no We have the following condition.

\bs\no {\bf 4.-} For each circuit $b_0b_1\cdots b_n$ in $B_N$, choose 
a two--block $b_ib_{i+1}$. The choice in {\bf 3} has to be such that
for each $a\in E_{b_i}$, $D_{a,b_{i+1}}:=E_{b_{i+1}}$.

\no {\bf 5.-} Define the transition matrix $M:A\times A\to\{0,1\}$ 
such that
%%%%%%%%%%%%%%%%%%%%%%%%%%%%%%%%%%%%%%%%%%%%%%%%%%%%%%%%%%%%%%%%%%%
\begin{equation}
M(a,a')=\left\{\begin{array}{lc}
 1 & \hbox{ if } a\in E_{b},\ a'\in D_{a,b'}\hbox{ and } N(b,b')=1,\\
 0 & \hbox{ otherwise. } \end{array}\right. \nonumber
\end{equation}
%%%%%%%%%%%%%%%%%%%%%%%%%%%%%%%%%%%%%%%%%%%%%%%%%%%%%%%%%%%%%%%%%%

\bs\no The factor map $\pi:A\to B$ such that $\pi^{-1}(b)=E_b$ for
each $b\in B$ has all the desired properties.

\bs\no For the class of examples constructed in this way, any
Markov measure in $A_M$ with support in the whole symbolic
set, induces a Gibbsian measure in the factor system $B_\pi := B_N$.

\bs \no One concrete instance of this kind is the topological
Markov chain $\{1,\ldots,5\}_M$, described by the digraph
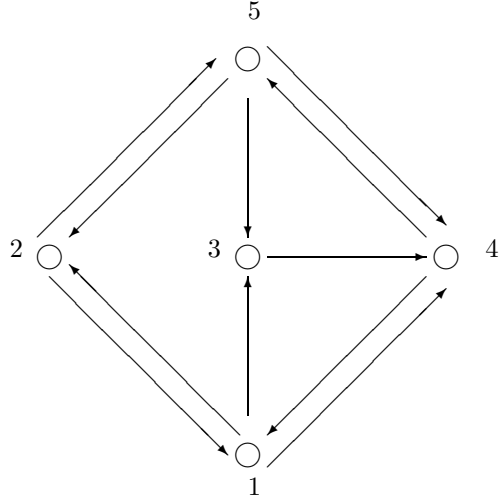
\begin{center}
\unitlength 1.5pt
\begin{picture}(200,130)(0,40)
\put(100,50){\circle{6}} \put(100,40){$1$}
\put(105,47){\vector(1,1){45}} \put(98,55){\vector(-1,1){43}}
\put(100,60){\vector(0,1){35}}
\put(50,100){\circle{6}} \put(40,100){$2$}
\put(100,100){\circle{6}}\put(90,100){$3$}
\put(100,150){\circle{6}}\put(100,160){$5$}
\put(50,95){\vector(1,-1){45}} \put(47,105){\vector(1,1){45}}
\put(105,100){\vector(1,0){40}}
\put(145,95){\vector(-1,-1){40}} \put(145,105){\vector(-1,1){40}}
\put(95,145){\vector(-1,-1){40}} \put(105,153){\vector(1,-1){45}}
\put(100,140){\vector(0,-1){35}}
\put(150,100){\circle{6}}\put(160,100){$4$}
\end{picture}
\end{center}

\bs \no Together with $\{1,\ldots,5\}_M$, consider the factor map
$\pi:\{1,\ldots,5\}\to\{a,b,c\}$, such that $\pi^{-1}(a)=\{1,5\}$,
$\pi^{-1}(b)=\{2,4\}$, and $\pi^{-1}(c)=3$.
The factor system $\{a,b,c\}_\pi$ is described by the digraph
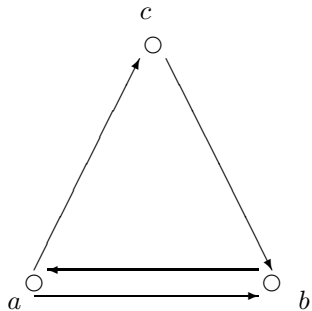
\begin{center}
\begin{picture}(200,150)(-50,0)
\put(10,10){\circle{6}} \put(0,0){$a$}
\put(100,10){\circle{6}}\put(110,0){$b$} \put(55,100){\circle{6}}
\put(50,110){$c$}
\put(10,5){\vector(1,0){85}} \put(95,15){\vector(-1,0){80}}
\put(10,15){\vector(1,2){40}}
\put(60,95){\vector(1,-2){40}}
\end{picture}
\end{center}

\bs \no Notice that one has only two periodic points with period $\leq 3$,
namely $(ab)^\infty$ and $(ba)^{\infty}$,
and for any Markovian measure with support in
$\{1,\ldots,5\}_M$, the corresponding matrices $\mm_{ab}$ and $\mm_{ba}$ 
are both $>0$.

%%%%%%%%%%%%%%%%%%%%%%%%%%%%%%%%%%%%%%%%%%%%%%%%%%%%%%%%%%%%%%%%%%%%%%%%%%%%%%%%%%%%%%%%%%%%%%%%%%%%%%%%%%%%%%%%%%
\subsection{The induced potential is not of finite range even when the original TMC is a full shift.} \label{fullshift}

\no The purpose of this section is to see what happens when there are no forbidden
blocks in the original system, i.e. it is a full shift. It is obvious that the
factor system is also a full shift. It will turn out from the example considered hereafter
that the image measure of a Markov measure supported by an arbitrary full shift is always a BGM.
We are rather interested in the possibility that the induced potential be of finite range. Our
example shows that it is indeed possible but somewhat exceptional.

\no Let $\mu_\phi$ be a Markovian measure for the full shift
$\{a,b,c,d\}^{\nn}$, and consider the mapping $\pi:\{a,b,c,d\}\to\{0,1\}$
such that  
%%%%%%%%%%%%%%%%%%%%%%%%%%%%%%%%%%%%%%%%%%%%%%%%%%%%%%%%%%%%%%%%%%
\begin{equation}
E_0:=\pi^{-1}(0)=\{a,b\}\  \text{ and }  E_1:=\pi^{-1}(1)=\{c,d\}.
\nonumber
\end{equation}
%%%%%%%%%%%%%%%%%%%%%%%%%%%%%%%%%%%%%%%%%%%%%%%%%%%%%%%%%%%%%%%%%

\no This mapping defines a factor map  
$\pi:\{a,b,c,d\}^{\nn}\to \{0,1\}^{\nn}$. 

\no The induced Gibbs measure has potential $\psi:\{0,1\}^{\nn}\to \rr$, 
which is completely determined by the $2\times 2$ positive matrices 
$\mm_{00}, \mm_{01}, \mm_{10}, \mm_{11}$, and by the $2\times 1$ vectors 
$\mu_0$ and $\mu_1$. 

\no Indeed, according to the Main theorem (whose hypothesis are
trivially satisfied !),
%%%%%%%%%%%%%%%%%%%%%%%%%%%%%%%%%%%%%%%%%%%%%%%%%%%%%%%%%%%%%%%%%
\begin{eqnarray}
\psi(\b)&=&\lim_{n\to\infty}\log(\1^{\dag}\mm_{b(0:1)}x_{b(1:n)})
\nonumber \\
        &=&\log(\1^{\dag}\mm_{b(0:1)}x_{b(1:\infty)}) \nonumber \\
        &=&\log(\1^{\dag}\mm_{b(0:1)}F_{b(1:k)}(x_{b(k:\infty)})),
           \ \forall k\in\nn.
\nonumber
\end{eqnarray}
%%%%%%%%%%%%%%%%%%%%%%%%%%%%%%%%%%%%%%%%%%%%%%%%%%%%%%%%%%%%%%%%%%

\no Let us remind that for all $1\leq m< n$,  
%%%%%%%%%%%%%%%%%%%%%%%%%%%%%%%%%%%%%%%%%%%%%%%%%%%%%%%%%%%%%%%%%
\begin{equation}
x_{b(m:n)}:=F_{b(m:n)}\left(\hat{\mu}_{b(n)}\right)=
F_{b(m:m+1)}\circ F_{b(m+1:m+2)}\circ \cdots 
\circ F_{b(n-1:n)}\left(\hat{\mu}_{b(n)}\right),
\nonumber
\end{equation}
%%%%%%%%%%%%%%%%%%%%%%%%%%%%%%%%%%%%%%%%%%%%%%%%%%%%%%%%%%%%%%%%%

\no with $\hat{\mu}_{\sb(n)}:=\mu_{\sb(n)}/|\mu_{\sb(n)}|_1$
%%%%%%%%%%%%%%%%%%%%%%%%%%%%%%%%%%%%%%%%%%%%%%%%%%%%%%%%%%%%%%%%
\begin{equation}
F_{ee'}:\Delta_{e'}\to\Delta_e, \ \text{ such that }\ 
F_{ee'}(x)=\frac{\mm_{ee'}x}{|\mm_{ee'}|_1}.
\nonumber
\end{equation}
%%%%%%%%%%%%%%%%%%%%%%%%%%%%%%%%%%%%%%%%%%%%%%%%%%%%%%%%%%%%%%%

\no In the case we treat now, all these transformations are
pure contractions, hence the limit 
%%%%%%%%%%%%%%%%%%%%%%%%%%%%%%%%%%%%%%%%%%%%%%%%%%%%%%%%%%%%%%%
\begin{equation}
x_{b(k:\infty)}:=\lim_{n\to\infty}F_{b(k:n)}
        \left(\hat{\mu}_{b(n)}\right)
\nonumber
\end{equation}
%%%%%%%%%%%%%%%%%%%%%%%%%%%%%%%%%%%%%%%%%%%%%%%%%%%%%%%%%%%%%%%%

\no exists for all $k\geq 1$. Notice also that the two simplices 
$\Delta_0$ and $\Delta_1$, are equivalent. Because of 
this, the functions $F_{ee'}$ can be considered as self--maps 
in the one--dimensional simplex
%%%%%%%%%%%%%%%%%%%%%%%%%%%%%%%%%%%%%%%%%%%%%%%%%%%%%%%%%%%%%%%%
\begin{equation}
\Delta:=\{x\in (0,1)\times (0,1):\ x(0)+x(1)=1\}.
\nonumber
\end{equation}
%%%%%%%%%%%%%%%%%%%%%%%%%%%%%%%%%%%%%%%%%%%%%%%%%%%%%%%%%%%

\no For each $e\in\{0,1\}$, the set of limit points 
%%%%%%%%%%%%%%%%%%%%%%%%%%%%%%%%%%%%%%%%%%%%%%%%%%%%%%%%%%%%%%
\begin{equation}\label{fractal}
{\mathcal F}_e:=\{x_{\se(1:\infty)}:\ \e\in \{0,1\}^{\nn}, \e(1)=e\},
\end{equation}
%%%%%%%%%%%%%%%%%%%%%%%%%%%%%%%%%%%%%%%%%%%%%%%%%%%%%%%%%%%%%

\no can be thought as a subset of the fractal limit of the Iterated 
System of Functions $(\Delta,\{F_{ee'}:\ ee'\in \{0,1\}\times\{0,1\}\})$. 

\no In general ${\mathcal F}_e$ is an uncountable set. For this it is 
enough that the fixed points of the mappings $F_{00}$ and $F_{11}$ be
different.  

\no We are in the situation where the values of $\psi$ can be obtained
through linear functionals on $\rr^2$, acting on the fractal set 
${\mathcal F}\subset \Delta\subset \rr^2$. 
For each couple $ee'\in \{0,1\}\times\{0,1\}$, these functionals are
defined by $x\mapsto \1^{\dag}\mm_{ee'}x$.

\bs
%%%%%%%%%%%%%%%%%%%%%%%%%%%%%%%%%%%%%%%%%%%%%%%%%%%%%%%%%%%%%%%%
\begin{proposition}\label{projectionofafullshift}
If the induced potential $\psi:\{0,1\}^{\nn}$ is of finite range, then at 
least one of the following conditions must hold.
\begin{itemize}
\item[(1)] The matrices $\mm_{00}$ and $\mm_{11}$ have the same positive
eigenvector.
\item[(2)] One of the matrices $\mm_{00}$, $\mm_{01}$, $\mm_{10}$, or 
$\mm_{11}$, is of rank 1.
\item[(3)] The vector $\1^{\dag}$ of dimension $1\times 2$ is a
left eigenvector for each one of the matrices 
$\mm_{00}$,$\mm_{01}$,$\mm_{10}$, and $\mm_{11}$.
\end{itemize} 
\end{proposition}
%%%%%%%%%%%%%%%%%%%%%%%%%%%%%%%%%%%%%%%%%%%%%%%%%%%%%%%%%%%%%

\bs\no Hence, if the Markov measure 
$\mu_\phi\in\bb\left(\{a,b,c,d\}^{\nn}\right)$ is such that,
under the factor map $\pi:\{a,b,c,d\}\to\{0,1\}$,
none of the hypotheses of the previous proposition holds, then the 
induced Gibbs measure $\nu_\psi\in\bb\left(\{0,1\}^{\nn}\right)$ 
cannot have a potential of finite range. The space of parameters 
defining a Markov measure $\mu_\phi\in\bb\left(\{a,b,c,d\}^{\nn}\right)$ 
is the Cartesian product of 4 simplices of dimension 3 (the 4 columns 
of the probability transition matrix). Conditions (1)-(3) 
of the previous proposition define a submanifold of dimension not 
greater that 9 inside that space of parameters: 1 dimension for 
condition (1), 4 dimensions for condition (2), and 4 dimensions for 
condition (3). Hence, a Markov measure $\mu_\phi\in\bb\left(\{a,b,c,d\}^{\nn}\right)$
generically induces a Gibbs measure whose potential
cannot have finite range.

%%%%%%%%%%%%%%%%%%%%%%%%%%%%%%%%%%%%%%%%%%%%%%%%%%%%%%%%%%%%%%%%%%%%%%%%%%%
\subsection{The induced potential at periodic points}\label{periodic}

\no It is worth to notice that the limit~(\ref{the-limit}) can
be effectively computed on periodic points.

\bs\no Let $\b\in B_\pi$ a periodic point of period $p\geq 1$. The set
%%%%%%%%%%%%%%%%%%%%%%%%%%%%%%%%%%%%%%%%%%%%%%%%%%%%%%%%%%%
\begin{equation}
\Omega_{\sb} := \{ \a\in A_M: \pi \a =\b\}
\nonumber
\end{equation}
%%%%%%%%%%%%%%%%%%%%%%%%%%%%%%%%%%%%%%%%%%%%%%%%%%%%%%%%%%%%%%
\no together with $\sigma^p$, define the full shift on $E_{\sb(0)}$.
Indeed, the preimage 
$\pi^{-1}(\b)$ of a periodic point $\b\in {\rm Per}_p(B_\pi)$, is a
$\sigma^p$-invariant subset of $A_M$. The system $(\pi^{-1}(\b),\sigma^p)$
is a TMC whose transition matrix is compatible with 
$M_{\sb(0:p)}:=\prod_{i=0}^{p-1}M_{\sb(i:i+1)}$.
By the hypothesis {\bf H2}, $M_{\sb(0:p)}$ is positive, and hence the
system $(\pi^{-1}(\b),\sigma^p)$ is a full 
shift on $E_{\sb(0)}:=\pi^{-1}(\b(0))$.

The topological pressure $P_{\Omega_{\sb}}(\phi;\sigma^p)$ of this
system, with respect to the potential $\phi$, is the logarithm of
the maximal eigenvalue of the matrix $\mm_{\sb(0:p)}$
(as defined in \eqref{Mmn}). Let us denote by
$G_{\sb(0:p)}$ and $D_{\sb(0:p)}$ the left and right eigenvectors 
associated with
$\rho:=\exp\left[P_{\Omega_{\sb}}(\phi;\sigma^p)\right]$.
They are chosen in order that  $G^{\dag}_{\sb(0:p)} D_{\sb(0:p)}=1$. 
Further set
$\hat{D}_{\sb(0:p)}:=D_{\sb(0:p)}/ \vert D_{\sb(0:p)}\vert_1$.
%%%%%%%%%%%%%%%%%%%%%%%%%%%%%%%%%%%%%%%%%%%%%%%%%%%%%%%%%%%%%%%%

\bs
%%%%%%%%%%%%%%%%%%%%%%%%%%%%%%%%%%%%%%%%%%%%%%%%%%%%%%%%%%%%%%%%%%%%%%%
\begin{proposition}\label{proj-period}
Let $\b\in B_\pi$ a periodic point of 
period $p\geq 1$ such that $\mm_{\sb(0:p)}$ is primitive, which is true, 
in particular, when hypothesis {\bf H2} is satisfied. Then
%%%%%%%%%%%%%%%%%%%%%%%%%%%%%%%%%%%%%%%%%%%%%%%%%%%%%%%%%%%%%%%%%
\begin{equation}\label{psi-period}
\psi(\b)= P_{\Omega_{\sb}}(\phi;\sigma^p)-
\log\left(\left\vert\mm_{\sb(1:p)}\hat{D}_{\sb(0:p)}\right\vert_1\right)\,.
\end{equation}
\end{proposition}
%%%%%%%%%%%%%%%%%%%%%%%%%%%%%%%%%%%%%%%%%%%%%%%%%%%%%%%%%%%%%%%

\no Proposition \ref{proj-period} leads to the following 
approximation formula. For any $\b'\in B_\pi$ and each $n\in \nn$ 
let $\b$ be a periodic point of minimal period $p(\b,n)$ in 
$[\b'(0:n)]$. Note that $p(\b,n)\to\infty$ as $n\to\infty$. 
It follows from \eqref{psi-is-holder} that 
$$
\vert \psi(\b')-\psi(\b)\vert \leq C\times \theta^{p(\b,n)},
$$
using for $\psi(\b)$ formula \eqref{psi-period}.

\no Let us recall that one can get a large class of Gibbs measure as a weak$^*$ limit, as $p\to\infty$,
of measures concentrated on $p$-periodic orbits obtained by counting
$p$-periodic orbits weighted by the potential, see e.g. \cite{KH}.

\no Another remarkable property of periodic orbits appears in Livsic's periodic points
theorem (see \cite{pp}): periodic points determine completely the
cohomology class of a H{\"o}lder continuous potential.

%%%%%%%%%%%%%%%%%%%%%%%%%%%%%%%%%%%%%%%%%%%%%%%%%%%%%%%%%%%%%%%%%%%%%%%%%%%%%
\subsection{Example of a non-gibbsian induced measure}\label{non-gibbs}

\no In this section we show that hypothesis {\bf H2} is essential to get a
well-defined potential on the {\bf whole} factor system. A concrete and simple example
is built such that at some point the induced potential is not defined. More precisely, this means
the sequence appearing in \eqref{psiasx} does not converge.
Let us emphasize that this shows the subtle effect produced by the presence of forbidden blocks in the
original system. Remember that without forbidden blocks one always gets a potential which is well-defined
everywhere (subsection \ref{fullshift}). 

\no Consider the TMC $\{a,b,c,d,e,f\}_M$ described by the digraph
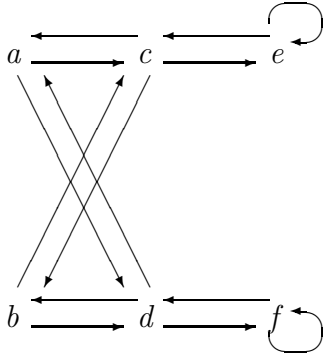
\begin{center} 
\begin{picture}(250,160)(50,20) 

\put(100,150){\large \em a} 
\put(100,50){\large \em b} 
\put(150,150){\large \em c} 
\put(150,50){\large \em d} 
\put(200,150){\large \em e} 
\put(200,50){\large \em f} 

\put(210,165){\oval(20,15)[tr]} 
\put(210,165){\oval(20,15)[tl]} 
\put(210,165){\oval(20,15)[br]} 
\put(208,158){\vector(-1,0){0}} 

\put(210,48){\oval(20,15)[bl]} 
\put(210,48){\oval(20,15)[br]} 
\put(210,48){\oval(20,15)[tr]} 
\put(208,56){\vector(-1,0){0}}  
 
\put(110,150){\vector(1,0){35}} 
\put(110,50){\vector(1,0){35}} 
\put(160,150){\vector(1,0){35}} 
\put(160,50){\vector(1,0){35}}

\put(150,160){\vector(-1,0){40}} 
\put(150,60){\vector(-1,0){40}} 
\put(200,160){\vector(-1,0){40}} 
\put(200,60){\vector(-1,0){40}} 

\put(105,145){\vector(1,-2){40}} 
\put(155,145){\vector(-1,-2){40}}
 
%\put(200,140){\vector(0,-1){70}}

\put(105,65){\vector(1,2){40}} 
\put(155,65){\vector(-1,2){40}} 

%\put(210,70){\vector(0,1){70}}

\end{picture} 
\end{center}

\noindent
The mapping
$\pi:\{a,b,c,d,e,f\}\to \{0,1\}$
such that $\pi^{-1}(0)=\{a,b,c,d\}$
and 
$\pi^{-1}(1)=\{e,f\}$ maps the TMC
$\{a,b,c,d,e,f\}_M$ onto the full shift $\{0,1\}^{\nn}$.

\bigskip
\noindent
Supply $\pi:\{a,b,c,d,e,f\}\to \{0,1\}$ with the with the 1-step Markov measure $\mu$, defined
by the probability transition matrix  
$$
\mm=\left(\begin{array}{cccccc}0     &  0    &  2\g  &   \g    &1-3\g  &   0   \\
                                0    &  0    &  \g   &   \g    &   0   & 1-2\g  \\
                                1/4  &  1/4  &  0    &   0    & 1/2   &   0   \\
                                1/4  &  1/4  &  0    &   0    &   0   & 1/2   \\
                                1/2  &  0    & 1-3\g  &   0   &3\g-1/2&   0   \\
                                0    & 1/2   &  0    & 1-\g   &  0    &2\g-1/2   
\end{array}\right),
$$
with $1/4 < \g < 1/3$. 
Since $\mm$ is double--stochastic, then the one--marginal is the uniform vector $(1/6)\1$.

\bigskip \noindent
Suppose that the induced measure $\nu=\mu\circ\pi^{-1}$ is a Gibbs measure defined by the 
potential $\psi:\{0,1\}\to\rr$. If this is the case, one must have
$\psi(0^\infty)=\limn\log\left(\1^{\dag}\mm_{00}x_{0^\infty(1:n)}\right)$,
with 
$$\mm_{00}=
           \left(\begin{array}{cccc}0    &  0    &  2\g  &  \g  \\
                                    0    &  0    &  \g   &  \g  \\
                                  1/4    & 1/4   &  0    &   0  \\
                                   1/4   & 1/4   &  0    &   0  \end{array}\right)
:=\left(\begin{array}{cc} {\bf 0}   &  A     \\
                              B    &{\bf 0} \end{array}\right)    
$$
and $x_{0^\infty(1:n)}=(\1^{\dag}\mm_{00}^{n-1}\1)^{-1}\mm_{00}^{n-1}\1$. 

\medskip \noindent
Since $\mm_{00}$ is a irreducible matrix whose second eigenvalue is zero, 
the Perron--Frobenius theory gives, for $k\geq 3$,
\begin{eqnarray}
\mm_{00}^{2k}=
\left(\begin{array}{cc}(AB)^k & {\bf 0}\\
                     {\bf 0}  & (BA)^k \end{array} \right) 
&=&\rho_{AB}^k 
\left(\begin{array}{cc}D_{AB}D_{BA}^{\dag} & {\bf 0}     \\
                     {\bf 0}  &D_{BA}D_{AB}^{\dag} \end{array} \right) \nonumber \\ 
\mm_{00}^{2k+1}=
\left(\begin{array}{cc}{\bf 0}  &(AB)^kA \\
                      (BA)^{k}B & {\bf 0}\end{array} \right)
\!\!\!&=&\!\!\!\rho_{AB}^k
\left(\begin{array}{cc}    {\bf 0}  & D_{AB}D_{BA}^{\dag}A \\
                 D_{BA}D_{AB}^{\dag }B & {\bf 0}       \end{array} \right),\nonumber 
\end{eqnarray}
where we used the notation previously used for maximal eigenvalues and associated eigenvectors
(preceding subsection). 
Notice that in this case $\rho_{AB}=\rho_{BA}=5\g/4$, $D_{AB}=\frac1{5}(3\ 2)^{\dag}$ and
$D_{BA}=(1\ 1)^{\dag}$.
With this we obtain
$$
\1^{\dag}\mm_{00}x_{{\bf 0}(1:2k+1)}=\frac{5\g+1}{4}, \quad 
\1^{\dag}\mm_{00}x_{{\bf 0}(1:2k+2)}=\frac{5\g}{5\g+1}.
$$
Since $(5\g+1)/4 \neq 5\g/(5\g+1)$ for $\g\neq 1/5$, then $\psi$ is not defined at the fixed 
point $0^\infty$. 
%Analogously, the same conclusion holds for any point of the form $y0^\infty$, 
%for any finite block $y\in \{0,1\}^k$, $k\geq 0$, 

%%%%%%%%%%%%%%%%%%%%%%%%%%%%%%%%%%%%%%%%%%%%%%%%%%%%%%%%%%%%%%%%%%%%
\subsection{The converse to Lemma \ref{row->markov} is false}
For this consider the topological Markov chain $A_M$, defined by the
digraph below, together with the mapping $\pi$ defined at the right
side of the picture. 
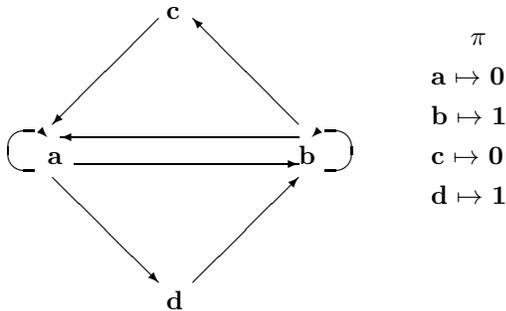
\begin{center}

\begin{picture}(250,200)(50,0)

\put(50,100){\oval(20,15)[bl]}
\put(50,100){\oval(20,15)[tl]}
\put(55,105){\vector(1,-1){0}}
\put(55,95){\bf a}
\put(160,100){\oval(20,15)[br]}
\put(160,100){\oval(20,15)[tr]}
\put(155,105){\vector(-1,-1){0}}
\put(150,95){\bf b}
\put(65,95){\vector(1,0){85}}
\put(150,105){\vector(-1,0){90}}
\put(150,110){\vector(-1,1){40}}
\put(100,150){\bf c}
\put(97,150){\vector(-1,-1){40}}
\put(57,90){\vector(1,-1){40}}
\put(100,40){\bf d}
\put(110,50){\vector(1,1){40}}  

\put(200,80){\hbox{{\bf d} $\mapsto$ {\bf 1}}}
\put(200,95){\hbox{{\bf c} $\mapsto$ {\bf 0}}}
\put(200,110){\hbox{{\bf b} $\mapsto$ {\bf 1}}}
\put(200,125){\hbox{{\bf a} $\mapsto$ {\bf 0}}}
\put(215,140){\hbox{$\pi$}}

\end{picture}

\end{center}

\no The factor system defined by $A_M$ and $\pi$
is the full shift $(\{0,1\}^{\nn},\sigma)$. 
For this note that the full shift $(\{a,b\}^{\nn},\sigma)$ is a subshift of $A_M$, 
and that the factor map restricted to this subshift is a conjugacy.

\no On the other hand, the submatrix $\mm_{01}$ is not a row allowable 
matrix since there is no symbol $e\in E_{1}(:=\pi^{-1}(1))$ such that $M(c,e)=1$. 
The same is true for $M_{10}$.  

%%%%%%%%%%%%%%%%%%%%%%%%%%%%%%%%%%%%%%%%%%%%%%%%%%%%%%%%%%%%%%%%
%%%%%%%%%%% NEW SECTION %%%%%%%%%%%%%%%%%%%%%%%%%%%%%%%%%%%%%%%%
\section{Concluding remarks and open questions}\label{open}

\no {\bf Some related works in ergodic theory.}
Some previous works \cite{boyle,KT} deal with the study of factor maps between TMC's
in the context of ergodic theory and dynamical systems. 
Let us mention the work of Walters~\cite{waltersbis},
where the concept of compensation function, which was first considered by
Boyle and Tuncel~\cite{boyle}, is used to characterize more general
factor maps that those considered in the present work.
Walters takes especially advantage of ideas and results from
the thermodynamic formalism of equilibrium measures. Let us also quote
the recent works \cite{shin} where the concept of compensation
functions is used to answer certain questions related to measures
maximizing some weighted entropy. We point out that in most cases all these
works study the behavior of Markov measures under {\it lifting} while in the
present work we were interested in {\it projecting} Markov measures.

\bs

\no {\bf Infinite-to-one factor maps and finite-to-one ones.}
There are two classes of factor maps between TMC's. Finite-to-one maps preserve the topological entropy
whereas infinite-to-one maps decrease it strictly. (An infinite-to-one map
is a map such that there is at least one point having an uncountable number of
preimages.) A simple combinatorial characterization
allows to determine whether a factor map is infinite-to-one: one has to check that
the map has a diamond. We refer the reader to \cite[Chap. 4]{K} for full details.
Hypothesis {\bf H2} in the Main theorem (section \ref{mainresult}) implies that $\pi$ is an infinite-to-one
factor map. A necessary condition to have an infinite-to-one factor map
is that $\#A>\#B$. But this not sufficient (see the nice example in
\cite[p. 97]{K}).

\bs

\no {\bf When the factor map is not a topological Markov map}.
In general a factor map maps a TMC to a strictly sofic subshift, see \cite{K}.
BGM's are in fact well-defined on sofic subshifts \cite{bowenbis}. We conjecture that ``nice''
factor maps should also map Markov measures to BGM's. We were forced to consider a subclass of
topological Markov factor maps, namely full row allowable factor maps (hypothesis {\bf H1} in the Main
theorem (section \ref{mainresult})). We believe that this is not a necessary condition but we are not
able to prove anything by using our present tools.
 
\bs

\no {\bf About rational probability measures and semi-group measures}.
In~\cite{hansel}, the authors introduce the concept of rational 
probability measures which are characterized in terms of formal power series.
In particular, they show that these measures are exactly the
measures obtained by the action of 1-block factor maps (alphabetic monoid morphisms
in their context) on 1-step Markov chains. Therefore the
Main theorem gives some sufficient conditions for a rational
probability measure to be a BGM. The same could be said for semi-group measures
that were introduced in \cite{KT}. Moreover the example of section \ref{non-gibbs}
shows that there are rational probablity measures or semi-group measures that are
not Gibbs measures.

\bs

\no {\bf Grouped Markov chains}. Our main result generalizes a result by Harris \cite{harris}
in the context of chains of infinite order (or chains with complete connections). This
author calls a grouped Markov chain what we call a projected Markov measure
and only considers full shifts, that is, there are no forbidden blocks.
Rephrased in his language, we can say that under our
hypothesis a grouped Markov shift is continuous with respect to its past with an
exponential continuity rate.

\bs

\no {\bf Equilibrium measures with a non-H{\"o}lder potential}. One can relax the hypothesis
of H{\"o}lder continuity of a potential function and still have a unique equilibrium state
satisfying property \eqref{BGI}. This is the case when the variation of the potential on cylinders
is not exponential (as in the H{\"o}lder case) but, for instance, summable.
A glance at our proof
shows that the H{\"o}lder continuity of the induced potential follows from the exponential
convergence in \eqref{the-limit}, see formula \eqref{potential-as-limit}. This
shows that by using the projective distance we can only obtain a H{\"o}lder continuous
induced potential.

\bs

\no {\bf Weak Gibbs measures and hypothesis {\bf H2}}. The example in section \ref{non-gibbs}
provides a simple example of the non-existence of the induced potential at {\it some
point}. This happens because there is the matrix $\mm_{00}$ associated to the fixed
point $0^\infty$ which is not strictly positive, in violation with hypothesis {\bf H2}
of the Main theorem.
It may happen that if {\bf H2} is satisfied for at least one periodic point the image
measure $\nu$ could be a BGM on a subset of $B_\pi$ of full measure, that is, the Bowen-Gibbs
inequality \eqref{BGI} could hold for $\nu$-almost all $\b\in B_\pi$. This situation has
been studied, see e.g. \cite{MRTvMV,yuri}, and such a measure is referred to as a weak Gibbs
measure.

%%%%%%%%%%%%%%%%%%%%%%%%%%%%%%%%%%%%%%%%%%%%%%%%%%%%%%%%%%%%%%%%
%%%%%%%%%%%%%%%%%%%%%%%% NEW SECTION %%%%%%%%%%%%%%%%%%%%%%%%%%%
%%%%%%%%%%%%%%%%%%%%%%%%%%%%%%%%%%%%%%%%%%%%%%%%%%%%%%%%%%%%%%%%
\section{Proof of some lemmas and propositions}\label{some-lemmas}

\no {\bf Proof of Lemma \ref{factorizing-b}.}

\no Write $\b$ as the concatenation $\b:=\b(0:\#B)\b(\#B+1:(\#B+1)+\#B)\cdots$. 
Inside each one of the factor blocks $\b(k(\#B+1):k(\#B+1)+\#B)$ there 
is at least one symbol appearing twice, i.~e., there are integers
$k(\#B+1)\leq m(k) < \ell(k)\leq k(\#B+1)+\#B$ such that
$\b(m(k))=\b(\ell(k))$. Hence, the lemma is proved.

%%%%%%%%%%%%%%%%%%%%%%%%%%%%%%%%%%%%%%%%%%%%%%%%%%%%%%%%%%%%%%%%
\bs

\no {\bf Proof of Lemma \ref{gibbsian-lemma}.}

\no From the hypothesis, for all
$\b\in \Omega$ and all $n\in\nn_0$ one gets
%%%%%%%%%%%%%%%%%%%%%%%%%%%%%%%%%%%%%%%%%%%%%%%%%%%%%%%%%
\begin{equation}
\exp\left(-\sum_{k=0}^nc(k)\right)\leq
\frac{\nu[\b(0:n)]}{\exp\left(\sum_{j=0}^{n}\psi(\sigma^j(\b))\right)}
\leq  \exp\left(\sum_{k=0}^n c(k)\right). \nonumber
\end{equation}
%%%%%%%%%%%%%%%%%%%%%%%%%%%%%%%%%%%%%%%%%%%%%%%%%%%%%%%%%

\no Since $\{c(n)\}_{n=0}^\infty$ is summable, then the Bowen-Gibbs
inequality~(\ref{BGI}) 
holds with a constant $K:=\sum_{n=0}^\infty c(n)$ and we are done.

%%%%%%%%%%%%%%%%%%%%%%%%%%%%%%%%%%%%%%%%%%%%%%%%%%%%%%%%%%%%%

\bs

\no {\bf Proof of Lemma \ref{row->markov}.}

\no Define the transition matrix $N:B\times B\to\{0,1\}$ such that
%%%%%%%%%%%%%%%%%%%%%%%%%%%%%%%%%%%%%%%%%%%%%%%%%%%%%%%%%%%%%%
\begin{equation}
N(b,b')=\left\{\begin{array}{lr} 
1 & \ \text{if} \ \exists\ a, a'\ M(a,a')=1, 
\ \pi(a)=b \ \text{and} \ \pi(a')=b',\\
0 & \ \text{otherwise}. \end{array} \right.
\nonumber
%%%%%%%%%%%%%%%%%%%%%%%%%%%%%%%%%%%%%%%%%%%%%%%%%%%%%%%%%%%%%%%
\end{equation}

\no We shall prove that the factor subshift $B_\pi$ and the TMC associated
to $B_N$ indeed coincide.

\bs\no
If $\b\in B_\pi$, then there exists $\a\in A_M$ such that 
$\pi(\a)=\b$. In particular, for each $i\in \nn_0$ the block 
$\a(i,i+1)$ satisfies 
%%%%%%%%%%%%%%%%%%%%%%%%%%%%%%%%%%%%%%%%%%%%%%%%%%%%%%%%%%%%%%%%%
\begin{equation}
M(\a(i),\a(i+1))=1\  \text{ and }\ \pi(\a(i))=\b(i) \ \forall i\in \nn_0.
\nonumber
\end{equation}
%%%%%%%%%%%%%%%%%%%%%%%%%%%%%%%%%%%%%%%%%%%%%%%%%%%%%%%%%%%%%%%%%%

\no Thus, $N(\b(i),\b(i+1))=1$ for each $i\in \nn_0$, and hence 
$\b\in B_N$. 

\bs\no On the other hand, if $\b\in B_N$, then for each $i\in \nn_0$ there
exists a block $a_ia_i'$ such that 
%%%%%%%%%%%%%%%%%%%%%%%%%%%%%%%%%%%%%%%%%%%%%%%%%%%%%%%%%%%%%%%%%%
\begin{equation}
M(a_i,a_i')=1,\ \  \pi(a_i)=\b(i) \ \text{ and }\ \pi(a_i')=\b(i+1).
\nonumber
\end{equation}
%%%%%%%%%%%%%%%%%%%%%%%%%%%%%%%%%%%%%%%%%%%%%%%%%%%%%%%%%%%%%%%%%%%

\no In general there is no reason that $a_i'=a_{i+1}$, but
since $M_{\sb(i:i+1)}$ is row allowable (recall Definition
\ref{row-allow}), given $a_i$ and $a_i'$, there exists 
$a_i''\in E_{\sb(i+2)}$ such that $M(a_i',a_i'')=1$. 
Then we can choose $a_{i+1}=a_i'$ and $a_{i+1}'=a_i''$. 
This choice is such that $a_ia_i'a_i''\equiv a_ia_{i+1}a_{i+1}'$
is $A_M$-admissible and $\pi(a_ia_{i+1}a_{i+1}')=\b(i:i+2)$. 
Thus, starting with $i=0$, we can proceed by induction in order to 
obtain a sequence $\a\in A_M$, such that $\pi(\a)=\b$ and $\a(i)=a_i$
for each $i\in \nn_0$. Therefore $\b\in B_\pi$. This concludes the proof.

%%%%%%%%%%%%%%%%%%%%%%%%%%%%%%%%%%%%%%%%%%%%%%%%%%%%%%%%%%%%%%%%%%%%%%%%%%%%%%%%%%%%%

\bs

\no {\bf Proof of Proposition \ref{projectionofafullshift}.}

\no Assume that the induced potential $\psi:\{0,1\}^{\nn}\to \rr$ is of
range $k$, for some $k\in \nn$. In that case, for all 
$\e\in\{0,1\}^{\nn}$, 
%%%%%%%%%%%%%%%%%%%%%%%%%%%%%%%%%%%%%%%%%%%%%%%%%%%%%%%%%%%%%%%
\begin{equation}
\exp(\psi(\e))=\1^{\dag}\mm_{(0:1)}F_{\se(1:k)}(x), \ \forall x\in
{\mathcal F}_{\se(k)}, 
\nonumber
\end{equation}
\no with ${\mathcal F}$ as defined in~\eqref{fractal}.

\no For this we have the following three logical possibilities: 
(i) either $\#{\mathcal F}_{\se(k)}=1$; or 
(ii) $\#{\mathcal F}_{\se(k)}>1$, but $F_{\se(1:k)}$ maps all points
in ${\mathcal F}_{\se(k)}$ to the same image; 
or (iii) $\#{\mathcal F}_{\e(k)} > 1$, $F_{\se(1:k)}$ maps two different
points in ${\mathcal F}_{\e(k)}$ to two different images, but
the linear functional $x\mapsto \1^{\dag}\mm_{\e(0:1)}x$ maps those
different images to the same value.

\bs\no If (i) holds, then the fixed points of the mappings 
$F_{00}$ and $F_{11}$, which belong to ${\mathcal F}_{\se(k)}$,
have to coincide. This means that $\mm_{00}$ and $\mm_{11}$ 
have the same positive eigenvector, and condition (1) in the
statement follows.

\bs\no If (ii) holds, then $\mm_{\se(1:k)}$ is necessarily a rank one 
matrix, and for this one of the matrices $\mm_{00}$, $\mm_{01}$,
$\mm_{10}$, or $\mm_{11}$, has to be of rank one. 
In this way condition (2) in the statement follows.

\bs\no Finally, if (iii) holds, we need 
$(\1^{\dag}\mm_{\se(0:1)})^{\dag}$ to be orthogonal to the simplex. 
In this case we have $\mm_{\se(0:1)}^{\dag}\1=\alpha \1$. 
Since $\e(0:1)$ is arbitrary, condition (1) follows. The proof is finished.

%%%%%%%%%%%%%%%%%%%%%%%%%%%%%%%%%%%%%%%%%%%%%%%%%%%%%%%%%%%%%%%%%%%%%%%%%%%%%%%%%%%%%

\bs

\no {\bf Proof of Proposition \ref{proj-period}.}

\no Since $\b$ is a periodic point of period $p$, \eqref{the-limit} becomes
%%%%%%%%%%%%%%%%%%%%%%%%%%%%%%%%%%%%%%%%%%%%%%%%%%%%%%%%%%%%%%%%%
\begin{equation}
\psi(\b) =\lim_{n\to \infty} 
\log\left(\frac{
\1^{\dag} \left( \mm_{\sb(0:p)}\right)^{\lfloor\frac{n}{p}\rfloor}
\mu_{\sb(0:n\ \textup{mod}\ p)} 
}{
\1^{\dag} \mm_{\sb(1:p)}\left( \mm_{\sb(0:p)}\right)^{\lfloor\frac{n}{p}\rfloor-1}
\mu_{\sb(0:n\ \textup{mod}\ p)}
}\right).
\nonumber
\end{equation}

\no To ease notation, let, for any $\b$ and $j\in\nn$,
$z_{\sb(0:j)}:= \mm_{\sb(0:j)}\mu_{\sb(j)}$.
Now apply Perron-Frobenius theorem \cite{seneta} to get
%%%%%%%%%%%%%%%%%%%%%%%%%%%%%%%%%%%%%%%%%%%%%%%%%%%%%%%%%%%%%%%%%%%%
$$
\psi(\b) =
$$
$$
\lim_{n\to \infty} \log\left(
\frac{
\1^{\dag} \rho^{\lfloor\frac{n}{p}\rfloor}
D_{\sb(0:p)} G^{\dag}_{\sb(0:p)} z_{\sb(0:n\ \textup{mod}\ p)} +
O\left(\lambda^{\lfloor\frac{n}{p}\rfloor}\right)
}{
\1^{\dag} \mm_{\sb(1:p)} \rho^{\lfloor\frac{n}{p}\rfloor-1}
D_{\sb(0:p)} G^{\dag}_{\sb(0:p)} z_{\sb(0:n\ \textup{mod}\ p)} +
O\left(\lambda^{\lfloor\frac{n}{p}\rfloor-1}\right)
}\right),
$$
%%%%%%%%%%%%%%%%%%%%%%%%%%%%%%%%%%%%%%%%%%%%%%%%%%%%%%%%%%%%%%%%%%%

\no where $\lambda$ is any number in $(\vert \lambda_2\vert, \rho)$
($\lambda_2$ is the eigenvalue of the next largest modulus after $\rho$).
Therefore
%%%%%%%%%%%%%%%%%%%%%%%%%%%%%%%%%%%%%%%%%%%%%%%%%%%%%%%%%%%%%%%%%%%
\begin{equation}
\psi(\b) = \log\rho +\lim_{n\to \infty} \log\left(
\frac{
\vert D_{\sb(0:p)}\vert_1
+O\left((\lambda/\rho)^{\lfloor\frac{n}{p}\rfloor}\right)
}{
\vert \mm_{\sb(1:p)} D_{\sb(0:p)}\vert_1 +
O\left((\lambda/\rho)^{\lfloor\frac{n}{p}\rfloor-1}\right)
}\right),
\nonumber
\end{equation}
%%%%%%%%%%%%%%%%%%%%%%%%%%%%%%%%%%%%%%%%%%%%%%%%%%%%%%%%%%%%%%%%%%%
and the proposition follows.
%%%%%%%%%%%%%%%%%%%%%%%%%%%%%%%%%%%%%%%%%%%%%%%%%%%%%%%%%%%%%%%%%

\bs

%%%%%%%%%%%%%%%%%%%%%%%%%%%%%%%%%%%%%%%%%%%%%%%%%%%%%%%%%%%%%%%%
%%%%%%%%%%%%%%%%%%%%%%%%%%%%%%%%%%%%%%%%%%%%%%%%%%%%%%%%%%%%%%%%
\no {\bf{\scriptsize ACKNOWLEDGMENTS}}.
We thank Karl Petersen and Fran{\c c}ois Blanchard for providing
us relevant references of related works when the first named author met 
them at the ``{\it Workshop on Dynamics and Randomness}''
held at Santiago, Chile (December 11--15, 2000). We are also grateful 
to the ZiF project {\it The Sciences of Complexity} for kind support.
The second named author was supported by ECOS/Nord-ANUIES program 
"Dynamics of extended systems".
We acknowledge the referees for their careful reading of the manuscript
and the remarks they made which lead to an improvement of the presentation.

%%%%%%%%%%%%%%%%%%%%%%%%%%%%%%%%%%%%%%%%%%%%%%%%%%%%%%%%%%%%%%%
%%%%%%%%%%%%%%%%%%%%%%%%%%%%%%%%%%%%%%%%%%%%%%%%%%%%%%%%%%%%%%%

\end{document}